\documentclass[11pt]{article}
\usepackage{graphicx}
\usepackage{psfrag}
\usepackage{amsmath,amsfonts,amsthm,amssymb}
\usepackage{amsfonts}
\usepackage{amsmath}
\usepackage{amssymb}
\usepackage{latexsym}

\usepackage{amssymb}
\usepackage{latexsym}

\newtheorem{Theorem}{Theorem}[section]
\newtheorem{Corollary}[Theorem]{Corollary}
\newtheorem{Lemma}[Theorem]{Lemma}
\newtheorem{Proposition}[Theorem]{Proposition}
\newtheorem{Conjecture}[Theorem]{Conjecture}
\newtheorem{Example}[Theorem]{Example}

\newtheorem{Remark}[Theorem]{Remark}
\newtheorem{Definition}[Theorem]{Definition}
\newtheorem{Claim}[Theorem]{Claim}

\newtheorem{Conclusion}[Theorem]{Conclusion}

\def\qed{{\hspace{2mm}{\hfill \small $\Box$}}}
\def\demo{ {\bf Proof.} }

\def\R{{\mathbb R}}
\def\H{{\mathbb H}}
\def\Z{{\mathbb Z}}
\def\Q{{\mathbb Q}}

\def\S{{\mathbb S}}
\def\E{{\mathbb E}}

\def\SS{{\mathcal S}}
\def\DD{{\mathcal D}}

\def\PP{{\mathcal P}}
\def\TT{{\mathcal T}}

\def\OO{{\mathcal O}}

\def\GG{{\mathcal G}}
\def\MM{{\mathcal M}}
\def\KK{{\mathcal K}}

\def\HH{{\mathcal H}}

\def\R{{\mathbb R}}
\def\H{{\mathbb H}}
\def\Z{{\mathbb Z}}
\def\Q{{\mathbb Q}}

\def\S{{\mathbb S}}
\def\E{{\mathbb E}}

\begin{document}
\title{Finiteness of $3$-manifolds associated with \\
non-zero degree mappings}

\author{Michel Boileau, J. Hyam Rubinstein and Shicheng
Wang\footnote{{\bf Acknowledgements} We thank a lot the referee for his carefull reading of our paper, for providing
many fine suggestions, and in particular for pointing out
the incompleteness of our induction argument in Section 5. We thanks also
 Y. Liu and H. Sun for their comments on our
drafts. The first and the last authors would like to acknowledge
support by MOST of China, and Institut Universitaire de France; and the last two
authors would like to acknowledge support by the Australian Research
Council}}
\date{}
\maketitle

\begin{abstract}
We prove a finiteness result for  the $\partial$-patterned guts
decomposition of all $3$-manifolds obtained by splitting a given
orientable, irreducible and $\partial$-irreducible 3-manifold along
a closed incompressible surface. Then using the Thurston norm, we
deduce that the JSJ-pieces of all 3-manifolds dominated by a given
compact 3-manifold belong, up to homeomorphism, to a finite
collection of compact 3-manifolds. We show also that any closed
orientable 3-manifold dominates only finitely many integral homology
spheres and any compact 3-manifolds orientable 3-manifold dominates
only finitely many exterior of knots in $S^3$.
\end{abstract}

\section{Introduction}

Maps between 3-manifolds have been studied for a long time, and
have become an especially active subject after Thurston's
revolutionary work on 3-manifold theory. The existence of non-zero
degree proper maps between compact orientable 3-manifolds is a
fundamental and difficult question in this area. We say that a
compact, orientable 3-manifold $M$ \emph{dominates} a compact
orientable 3-manifold $N$ if there is a non-zero degree proper map
$f: M \to N$. When the degree of $f$ is one, we say that $M$
\emph{1-dominates} $N$.

The following simple and natural question was raised in the 1980's
(and formally appeared in the 1990's, see (\cite[Problem 3.100
(Y.Rong)] {Ki}, and also \cite{W2}.

\noindent {\bf Question 1.} \emph{Does a closed orientable
$3$-manifold 1-dominate at most finitely many closed, irreducible
and  orientable $3$-manifolds?}

If we allow any degree, 3-manifolds supporting one of the
geometries $\S^3$, $\widetilde{PSL_2 (\R)}, Nil$  can dominate
infinitely many 3-manifolds. Thus any closed orientable
3-manifold which dominates such closed 3-manifolds indeed
dominates infinitely many 3-manifolds. At the moment these are the
only known examples, so the following generalization of Question 1
makes sense:

\noindent {\bf Question 2.} \emph{Let $M$ be a closed orientable
$3$-manifold. Does $M$ dominate at most finitely many closed,
irreducible, orientable $3$-manifolds $N$ not supporting the
geometries of $\S^3$, $\widetilde{PSL_2 (\R)}, Nil.$?}

Many related partial answers to Questions 1 and 2 have already
appeared in the literature, see for example \cite{Ro1}, \cite{Ro2},
\cite{Ro3}, \cite{BW}, \cite{HWZ1}, \cite{HWZ2}, \cite{RW},
\cite{So1}, \cite{So2}, \cite{So3}, \cite{Re}, \cite{WZ}, \cite{De},
\cite{De1}, \cite{Gu}, \cite{BCG}, \cite{BBW}.

Finiteness of closed irreducible targets implies
finiteness of possibly non prime closed targets: the number of prime factors in the target is bounded by the number of closed, disjoint, non parallel, essential surfaces in the domain; furthermore a
connected sum of finitely many closed 3-manifolds 1-dominates each
of its prime summands, which are either irreducible or homeomorphic to $S^1 \times S^2$.

Since a non-zero degree proper map between two compact orientable 3-manifolds induces a non-zero degree map between their doubles,  finiteness of closed
targets implies finiteness of the targets in the setting of compact
orientable 3-manifolds (see Remark \ref{compact} at the end of Section \ref{JSJ-pieces}).

First we introduce some standard terminology in 3-manifold topology,
see \cite{Ja}.

In this paper, all surfaces and 3-manifolds are compact and
orientable. Also we will work in the piecewise linear category, so
all spaces and maps will be PL. Suppose $S$ (resp. $P$) is a
properly embedded surface (resp. an embedded 3-manifold) in a
3-manifold $M$.  We use $M\setminus S$ (resp. $M\setminus P$) to
denote the resulting manifold obtained by splitting $M$ along $S$
(resp. removing $\text{int} P$, the interior of $P$). Note that
we allow the possibility that
$S$ is not connected, so that it has finitely many components.

A 3-manifold $M$ is:
\begin{itemize}
\item \emph{prime} if it is not the connected sum of two
3-manifolds neither of which is $S^3.$

\item \emph{irreducible} if every embedded sphere in $M$ bounds
a ball in $M$. A prime orientable 3-manifold which is not irreducible
is homeomorphic to $S^2 \times S^1$.

\item \emph{$\partial$-irreducible} if for every properly embedded
disc $D$ in $M$, there is a ball $B\subset M$ and a disc
$D'\subset\partial M$, such that $\partial D =
\partial D'$ and $\partial B=D\cup D'.$

\item \emph{atoroidal} if every $\Z \oplus \Z$ subgroup in
$\pi_1M$ is conjugate into $\pi_1\partial M$ and in addition
$\pi_1M$ is not virtually abelian. An irreducible orientable $3$-manifold such that every $\Z
\oplus \Z$ subgroup of $\pi_1M$ is conjugate into $\pi_1\partial
M$ is either atoroi\textit{}dal,  $T^2 \times (0,1)$, or the twisted I-bundle over the Klein
bottle.
\end{itemize}

The JSJ-decomposition (\cite{JS}, \cite{Joh}) of a compact orientable irreducible
3-manifold $M$ is the canonical splitting of $M$ along a finite
(possibly empty) collection $\TT$ of disjoint and non-parallel,
nor boundary-parallel, incompressible, embedded tori into maximal
Seifert fibered or atoroidal compact submanifolds.  We call the
components of $M\setminus \TT$ the \emph{JSJ-pieces} of $M$.

Thurston's geometrization conjecture stated that the atoroidal
JSJ-pieces support a hyperbolic or a spherical metric on their
interiors. W. Thurston proved his conjecture for any Haken
3-manifold (i.e. a compact, orientable, irreducible and
$\partial$-irreducible $3$-manifold which contains a properly
embedded essential surface; for details see \cite{Th1}, \cite{Th2}).
The full conjecture has been settled recently by G. Perelman
\cite{Per} (see \cite{KL}, \cite{MT}, \cite{CZ}, \cite{BBBMP}),
hence every compact orientable, irreducible 3-manifold is
\emph{geometrizable} in the sense that it satisfies Thurston's
geometrization conjecture.

The fact that the targets are geometrizable 3-manifolds is crucial
when considering Questions 1 and 2. A consequence is that a positive
answer to Question 2 implies a positive answer to Question 1.
Furthermore, since the results of \cite{So1} (see also \cite{Re},
\cite{Gu}, \cite{BCG}) and \cite{WZ} show that a closed orientable
3-manifold dominates only finitely many 3-manifolds supporting
either a hyperbolic structure with finite volume or a Seifert
geometry $\H^2\times \E^1$, Question 2 reduces to the following:

\noindent {\bf Question 3.} \emph{Let $M$ be a closed orientable
$3$-manifold. Does $M$ dominate at most finitely many, closed,
orientable, irreducible 3-manifolds $N$ with non-trivial JSJ
decomposition?}

There are some partial results for Question 3 in the case of
sequences of degree 1 maps (see \cite{Ro1}, \cite{So2}), or when the
domain and the target have the same simplicial volume (see
\cite{So3}, \cite{De}). Question 3 is solved in \cite{De} when $M$
is a graph manifold.

A general approach to Question 3 can be divided into two
steps:\begin{enumerate}

\item \emph{Finiteness of JSJ-pieces}: show that there is a finite
set $\HH\SS(M)$ of complete hyperbolic 3-manifolds with finite
volume and of Seifert manifolds such that each JSJ-piece of a
3-manifold $N$ dominated by $M$ belongs to $\HH\SS(M)$.

\item \emph{Finiteness of gluing}: For a given finite set
$\HH\SS(M)$ of complete hyperbolic 3-manifolds with finite volume
and of Seifert manifolds, there are only finitely many ways of
gluing elements in $\HH\SS(M)$ to get closed 3-manifolds dominated
by $M$.
\end{enumerate}

Notice that with our terminology, a manifold supporting a Sol
geometry has a non-trivial JSJ-decomposition with only one piece
homeomorphic to a product $T^2 \times I$ or two pieces homeomorphic
to the twisted I-bundle over the Klein bottle. For such manifolds in
the target, the finiteness of JSJ-pieces is trivially true, while
the finiteness of gluing is much more subtle (see \cite{WZ} for
1-domination  and \cite{BBW} if we allow arbitrary non-zero degree).

 T. Soma proved the finiteness of
hyperbolic JSJ-pieces in \cite{So2}. One of the main results of this
paper is to complete the proof of the first step by proving the
finiteness of the Seifert fibered JSJ-pieces:

\begin{Theorem}[Finiteness of JSJ-pieces] \label{local domination}
Let $M$ be a closed, orientable, 3-manifold. Then there is a
finite set $\HH\SS(M)$ of complete hyperbolic 3-manifolds with
finite volume and of Seifert fibered 3-manifolds, such that the
JSJ-pieces of any closed, orientable, irreducible
3-manifold $N$ dominated by $M$ belong to $\HH\SS(M)$, provided
that  $N$ does not support the geometries of $\S^3$,
$\widetilde{PSL_2 (\R)}, Nil.$
\end{Theorem}

The finiteness of the Seifert fibered JSJ-pieces follows from a
finiteness result for the Thurston norm of all compact 3-manifolds
$M_{S} = M\setminus S$, where $S$ runs over all incompressible,
orientable surfaces (not necessary connected) in $M$, see Section
\ref{finite norm}. This latter result is derived from the
finiteness of ``patterned guts" of all the
manifolds $M_{S} = M\setminus S$, where $S$ runs over all
incompressible, orientable surfaces in $M$, which we prove in
Section \ref{guts}.

We also prove the finiteness of gluing when the targets are
irreducible, integral homology 3-spheres. Together
with Theorem 1.1 (Theorem  \ref{local domination}) this gives a
positive answer to Question 3 when the targets are integral
homology spheres.

\begin{Theorem}\label{homology spheres} Any closed orientable $3$-manifold
dominates only finitely many integral homology
$3$-spheres.
\end{Theorem}

Without any restriction on the possible degrees of the maps or on
the geometry of the target, this is the best result one can
expect, since any closed orientable 3-manifold dominates all
3-dimensional lens spaces, which are rational homology spheres.

Since a degree-one map induces an epimorphism at the level of homology groups,
Theorem \ref{homology spheres} gives a positive answer to  Question 1 when the domain is an integral homology sphere:

\begin{Corollary}\label{integral spheres} An integral homology 3-sphere 1-dominates at most
finitely many closed 3-manifolds.
\end{Corollary}

The argument for integral homology spheres can be modified to prove
the following corollary.

\begin{Corollary}\label{knot spaces} Any compact orientable 3-manifold dominates at most
finitely many knot complements in $\S^3$.
\end{Corollary}

This corollary is related to a question of J. Simon on epimorphisms
between knot groups (see \cite[Problem 1.12 (J. Simon)]{Ki} and
Section \ref{knots}).

The paper is organized as follow: The finiteness of patterned guts
is discussed in Section \ref{guts}; the finiteness of Thurston norm
and Gromov volume is discussed in Section \ref{norm-volume}; the
finiteness of JSJ-pieces is proved in Section \ref{JSJ-pieces}. The
last Sections \ref{spheres} and \ref{knots}, are devoted to finite
domination results when the targets are integral homology 3-spheres
or knot complements in $\S^3$.

We end the introduction by the following

\begin{Remark} We could define the notion of domination between 3-manifolds
which are not necessarily orientable in terms of geometric degree
\cite{Ep}. But, then there are examples of non-orientable
(hyperbolic) 3-manifolds which 1-dominate infinitely many orientable
(hyperbolic) 3-manifolds (see \cite{Ro3}, \cite{BW}). In
\cite[Section 3]{BW}, , by lifting the maps in those examples to the
orientable double cover of the domain, maps between orientable
hyperbolic 3-manifolds are produced whose degree is $1+(-1)= 0$
rather than 2, as wrongly claimed there. This error has been pointed
out by T. Soma and many others.
\end{Remark}

\section{Finiteness of patterned guts}\label{guts}

In the Jaco-Shalen-Johannson decomposition of a compact orientable
3-manifold along essential tori and annuli, the guts consist of the
pieces which are not $I$-bundles over surfaces with negative Euler
characteristic. Finiteness of guts is a basic principle, which
originated from H. Kneser's work, see  for example \cite{A},
\cite{Ga2}, \cite{JR} for some recent applications of  guts in
3-manifold theory. We first introduce the notion of patterned guts
needed for our study of non-zero degree maps.

Suppose $X$ is a $\partial$-irreducible and irreducible, compact,
orientable 3-manifold. According to Jaco-Shalen-Johannson theory
([Ja], [JS], [Joh]), there is a unique, up to proper isotopy, characteristic 3-submanifold $\Sigma \subset X$ which is an union of Seifert spaces and $I$-bundles.

This characteristic submanifold has a unique decomposition, up to
proper isotopy:

$$\Sigma= (\Sigma\setminus IB^-_{X})\cup
IB^-_{X},$$

\noindent where $IB^-_{X}$ is formed by the components of the
Seifert pairs which are $I$-bundles over surfaces $F$, where $F$ has
negative Euler characteristic $\chi(F)$ if $\partial F\ne
\emptyset$. We make the following convention in this paper: if a
component of $X$ is a Seifert manifold and also an $I$-bundle over
a surface, we will always consider it as an $I$-bundle.

Therefore we have a decomposition

$$X = (X \setminus IB^-_{X})\cup_{A_{X}} IB^-_{X} = G_{X} \cup_{A_X} IB^-_{X},$$

\noindent where $A_{X}$ is the collection of frontier annuli of
$IB^-_{X}$ in $X$. We call $G_{X}=X \setminus IB^-_{X}$ the
\emph{guts} of $X$, and the decomposition above the
\emph{GI-decomposition} for $X$.  The embeddings of $G_{X}$,
$A_{X}$ and $IB^-_{X}$ are unique up to proper isotopy in $X$.

Suppose $S$ is a closed, incompressible surface in an irreducible
3-manifold $M$. Then $M_S=M\setminus S$ is $\partial$-irreducible
and irreducible. For such a surface $S$, we write the $GI$
decomposition of $M_S$ as

$$M_S = G_S \cup_ {A_S} IB^-_S.$$

\begin{Definition} Suppose $X$ is an orientable, irreducible and $\partial$-irreducible
3-manifold. A $\partial$-pattern for $X$ is a finite collection of
disjoint annuli $A \subset \partial X$, and given $A$ we say that
$X$ is $\partial$-patterned.
\end{Definition}

\begin{Example} For each component $G$ of $G_S$, $G\cap A_S$ is a
$\partial$-pattern for $G_S$. We often call the pair $(G, G\cap
A_S)$ a \emph{patterned guts component} for the surface $S$.
\end{Example}

The main result of this section is the following finiteness result
for patterned guts:

\begin{Theorem} [Patterned guts finiteness]\label{guts finiteness}
Let $M$ be a closed, orientable, irreducible 3-manifold. Then
there is a finite set $\GG(M)$ of connected, compact, orientable,
$\partial$-patterned 3-manifolds such that for each closed,
incompressible (not necessarily connected) surface $S\subset M$,
all patterned guts components of $(G_S, G_S\cap A_S)$ belong to
$\GG(M)$.
\end{Theorem}


\demo The proof of Theorem \ref{guts finiteness} consists of three
steps.

{\bf Step 1}. \emph{Construct a first ``approximation" to the
GI-decomposition by applying a refined Kneser argument.}

Fix a triangulation $K$ of $M$. Suppose that $K$ has $t$ tetrahedra.
For simplicity, we also assume that $K$ has only one vertex $v$ (see
\cite{JR} for example). Let $S_v$ be the normal sphere which is the
boundary of a small regular neighborhood $B_v$ of $v$. Suppose that
$S$ is a closed orientable incompressible surface in $M$. First
deform $S$ to be a normal surface in $(M,K)$. We can assume that
$S\cap S_v=\emptyset$. Let $S_*=S\cup S_v$.

Each tetrahedron $T$ has seven normal disc types, four triangular
types and three quadrilateral types, see Figure 1. Since $S_*$
contains $S_v$ and $S_*$ is embedded, for each tetrahedron $T$ of
$K$, $T\cap S_*$ contains all four triangular normal disc types
but at most one quadrilateral normal disc type.

 Let $M_*=M\setminus B_v$,  $K_*=K\cap M_*$, and
$T_*=T\cap M_*$ for each tetrahedron $T$ in $K$. Then $K_*$ is a
truncated triangulation of $M_*$, and each $T_*$ is a truncated
tetrahedron. Now we consider $S\subset M_*=|K_*|.$

If $S \cap T_*$ contains a quadrilateral normal disc, then
$T_*\setminus S$ contains two non-product regions, which are
truncated prisms : they are truncated from $T$  by using this
quadrilateral normal disc and four non-parallel triangular normal
discs $S \cap T$, see Figure 2. The boundary of each such a
truncated prism component has seven faces:

\noindent (1) two triangular normal discs  which lie in $S\cup
S_v$;

\noindent (2) one quadrilateral normal disc  which lies in $S$;

\noindent (3) two hexagonal faces which lie in the boundary of
$T$;

\noindent (4) two quadrilateral faces which lie in the boundary of
$T$.

\begin{center}%
\includegraphics[totalheight=8cm]{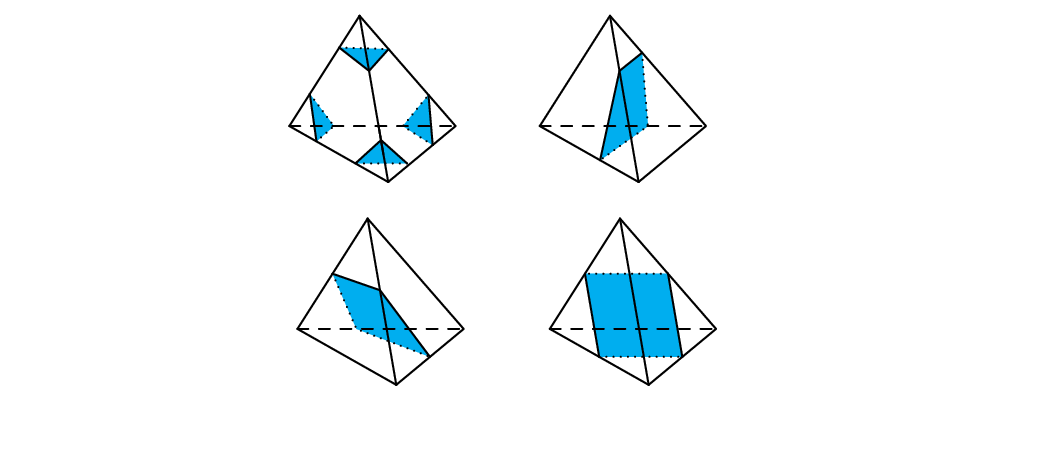}%

Figure 1
\end{center}
\begin{center}%
\includegraphics[totalheight=7cm]{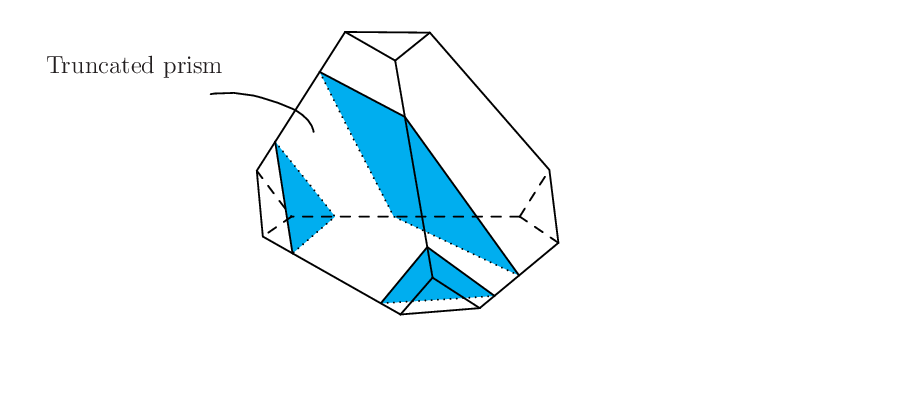}

Figure 2
\end{center}

If $S \cap T_*$ contains no quadrilateral normal disc, then
$T_*\setminus S$ contains just one non-product region, which is a
truncated tetrahedron: it is truncated from $T$  by using four
non-parallel normal discs of triangular type. The boundary of such
a truncated tetrahedron component has eight faces:

\noindent (5) four normal discs of triangular type which lie in
$S\cup S_v$;

\noindent (6) four hexagonal faces which lie in the boundary of
$T$.

Note that each remaining component of $T_*\setminus S$ is a
product region, whose boundary is formed by two normal discs of
the same triangular (resp. quadrilateral) type and three (resp.
four) vertical quadrilateral faces which lie in $\partial T$, see
Figure 3. Moreover in $K_*\setminus S$, each hexagonal face given
in (3) or (6) is identified with a hexagonal face given in (3) or
(6), and each quadrilateral face given  in (4) is either
identified with a quadrilateral face given  in (4), or with a
vertical quadrilateral face of a product region.

\begin{center}%
\includegraphics[totalheight=7cm]{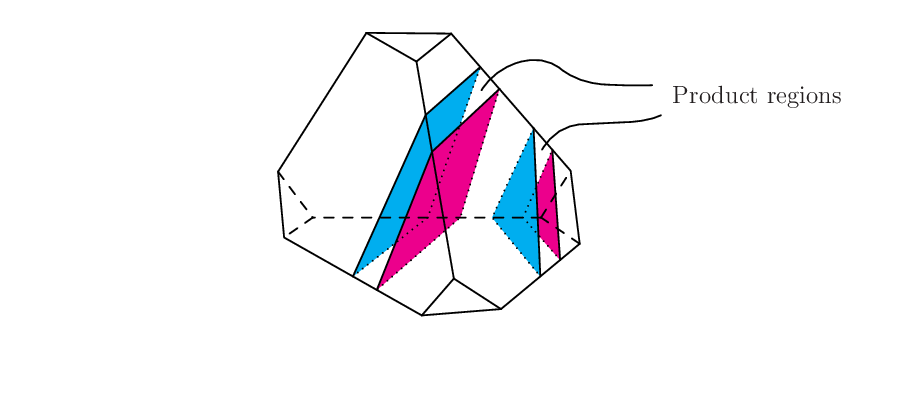}%

Figure 3
\end{center}

Let $Q$ be a quadrilateral face given in (4). If in $K_*\setminus
S$, $Q$ is identified with a vertical quadrilateral face of a
product region, we call $Q$ a \emph{frontier quadrilateral face}.
Otherwise we call $Q$ a \emph{non-frontier quadrilateral face}.

Now we glue together the truncated prism components and the
truncated tetrahedron components of $K_*\setminus S$ along their
hexagonal faces and their non-frontier quadrilateral faces to
get pieces $P_1$ ,...., $P_k$. Let $G^1_S$ be the union of those
pieces $P_i$, $i=1,...,k$.

Note that $\partial M_{*S}=S_1\cup S_2\cup S_v$, where $S_1$ and
$S_2$ are two copies of $S$, and
$$\partial G^1_S=(\partial
M_{*S}\cap G_S^1)\cup \text {(union of frontier
quadrilaterals).}\hskip 2 true cm (i)
$$

Now we have $M_*\setminus S =G^1_S\cup ((M_*\setminus S)\setminus
G^1_S)$. The components of $(M_*\setminus S)\setminus G^1_S$ are
obtained by gluing the product regions along their vertical
quadrilateral faces, hence they are $I$-bundles, whose union is
denoted by $IB^1_S$. The set $IB^1_S$ is a product or a twisted
$I$-bundle over a compact surface $S'$, denoted by $N(S')$. Let
$N(\partial S')$ denote the $I$-bundle structure restricted to
$\partial S'$. Then

$$\partial IB^{1}_S =(\partial
M_{*S}\cap \partial IB^{1}_{S}) \cup N(\partial S').\hskip
2 true cm (ii)$$

Clearly

$$\partial
M_{*_S} =(\partial M_{*S}\cap G_S^1)\cup (\partial M_{*S}\cap
\partial IB^{1}_{S}).\hskip
2 true cm (iii)$$

Combining the formulas (i) (ii) and (iii), it follows that the
annuli $N(\partial S')$ in (ii) are identified with the union of
frontier quadrilaterals in (i). In conclusion, all those frontier
quadrilaterals  form the intersection $IB^1_S\cap G^1_S$, which is
a union of finitely many properly embedded annuli in $M_*\setminus
S$, denoted by $A^1_S$. We call $A_S^1$ the frontier annuli of
$G_S^1$ (of $IB_S^1$). Now we get our first ``approximation"
$GI$-decomposition

$$M_*\setminus S=G^1_S\cup _{A^1_S} IB^1_S.$$

For each $S\subset M$, $G^1_S$ is constructed from $n \le t$
truncated tetrahedra and $m \le 2t$ truncated prisms by gluing
their hexagonal faces and non-frontier quadrilateral faces in
pairs. It follows that there is a bound for the combinatorial
(therefore the topological) types of the components of $G^1_S$. A
very crude bound is $5^t$, obtained by noting that there are $5$
choices for each tetrahedron consisting of the empty set, the
truncated tetrahedron or one of the 3 possible truncated prisms
(note if we have any quadrilateral type we always get two
truncated prisms in our guts). Moreover once $G_S^1$ is formed,
the position of the frontier annuli $A_S^1\subset G_S^1$ is fixed.
Hence we reach the following conclusion:

\begin{Conclusion}\label{conclusion1} There is a finite set
$\GG^1(M)$ of compact, orientable, connected, $\partial$-patterned
3-manifolds such that for each closed, incompressible surface
$S\subset M$, all patterned components of  $(G^1_S, G^1_S\cap
A_S^1)$ belong to $\GG^1(M)$.
\end{Conclusion}

\noindent {\bf Step 2.} \emph{Construct a second ``approximation" to
the GI-decomposition by absorbtion of ``tiny" patterned 3-manifolds
from $G^1_S\cup _{A^1_S} IB^1_S$.}

 Suppose that a component $A_i$ (or a pair of components $A_i$ and $A_j$) of $A^1_S$ separates
a component $P$ from $M_*\setminus S$ such that one of the following
patterned 3-manifolds occurs:

\noindent (i) $(P, A_i)=(D^2\times I, \partial D^2\times I)$, or
$((D^2\times I)\setminus B_v, \partial D^2\times I)$;

\noindent (ii) $(P, A_i)=(D^2\times S^1, I\times S^1)$  or $(P,
A_i)=((D^2\times S^1) \setminus B_v, I\times S^1)$, for some
interval $I \subset
\partial D^2$;

\noindent (iii) $(P, A_i\cup A_j)=(A\times I, \partial A\times I)$
or $((A\times I)\setminus B_v, \partial A\times I)$ for some
annulus $A$;

We call any patterned 3-manifold of one of the types above {\it
tiny}.

Note that a tiny patterned 3-manifold $P$  may contain other tiny
patterned 3-manifolds. Therefore $P$ may contain (finitely) many
components of $A_S^1$. But since $A_S^1$ has finitely many components,
there are only finitely many tiny patterned $P$.

Let $P$ be a tiny patterned 3-manifold. We eliminate $P$
by gluing it to its neighboring component(s) along $A_i$ (and
$A_j$) and then delete from $A_S^1$ all components of $A_S^1$
in $P$. In this manner, we also eliminate all tiny
patterned 3-manifolds contained in $P$. In such an absorbtion
process, we get a new decomposition $G_S^1(1)\cup
_{A_S^1(1)}IB_S^1(1)$ : all components in $G_S^1$ and $IB_S^1$
which are contained in $P$ or are adjacent to $P$ become a new component
of $G_S^1(1)\cup IB_S^1(1)$, all the remaining components in
$G_S^1$ and $IB_S^1$ are preserved, and ${A_S^1(1)}$ is obtained by
removing from $A_S^1$ all components of $A_S^1$ in $P$. The new
component in the new decomposition which contains $P$ is considered as a
``pseudo" $I$-bundle (respectively a ``pseudo" guts) component if and only
if the neighboring components of $P$ are $I$-bundles (respectively guts components).

Now consider the tiny patterned 3-manifolds of the decomposition
$G_S^1(1)\cup _{A_S^1(1)}IB_S^1(1)$ defined as above (which indeed
is a sub-collection of the tiny patterned 3-manifolds of
$G_S^1\cup _{A_S^1}IB_S^1$). If there are some, we can continue
this absorbtion process to get a new decomposition $G_S^1(2)\cup
_{A_S^1(2)}IB_S^1(2)$. Repeating this process we get a sequence of
decompositions $G_S^1(n)\cup _{A_S^1(n)}IB_S^1(n)$.  Since $A_S^1$
has  only finitely many components and the number of components of
$A_S^1(n)$ is strictly decreasing, this absorbtion process must
stop for some $n$. Then we get our second ``approximate"
$GI$-decomposition without tiny patterened 3-manifolds, which is
denoted by

$$M_*\setminus S=G^2_S\cup_{A^2_S} IB^2_S.$$

 Now we claim the following:

\begin{Claim} Each annulus in $A_S^2$ is incompressible and
$\partial$-incompressible in $M_*\setminus S$.
\end{Claim}

\demo Suppose some annulus $A_i \subset A_S^2$ is compressible in
$M_*\setminus S$. Since $M_*\setminus S$ is
$\partial$-irreducible, each component of $\partial A_i$ bounds a
disk in $\partial (M_*\setminus S)$. Since $M\setminus S$ is
irreducible, $A_i$ must separate from $M_*\setminus S$ either a
component homeomorphic to $D^2\times I$ or to $(D^2\times
I)\setminus B_v$. This contradicts the fact that no $A_i$ meets
the condition (i).

Suppose that some annulus $A_i \subset A_S^2$ is
$\partial$-compressible in $M_*\setminus S$. Since $M_*\setminus
S$ is irreducible and $\partial$-irreducible, it is not difficult
to verify that $A_i$ must separate from $M_*\setminus S$ a
component $P$ homeomorphic to a solid torus or a punctured solid
torus and which meets the condition (ii). This again gives a
contradiction.\qed

Since  $A_S^2$ is incompressible and $\partial$-incompressible in
$M_*\setminus S$,  $A^2_S$ does not meet $S_v$, and we can plug
the ball $B_v$ back into $M_*\setminus S$ to get a  new ``pseudo"
$GI$-decomposition for $M\setminus S$, still denoted as

$$M\setminus S=G^2_S\cup _{A^2_S} IB^2_S.$$

Let $m$ be the number of pattern annuli in $\GG^1(M)$. Let
$\PP(M)$ be the set of patterned 3-manifolds consisting of $m$
copies of a patterned 3-manifold of each type (i), (ii) and (iii),
and of one 3-ball. Then the patterned 3-manifolds obtained from
$\GG^1(M)$ and $\PP(M)$ by identifying some of their pattern
annuli in pairs, and possibly plugging in the 3-ball, is a finite set
$\GG^2(M)$ of patterned 3-manifolds. Since $G^2_S$ is obtained
from a subset of $G_S^1\subset\GG^1(M)$ and a subset of $\PP(M)$ by
identifying some of their pattern annuli in pairs, and possibly
plugging in the 3-ball, it follows that, up to homeomorphism, the
components of $G_S^2$ belong to $\GG^2(M)$. Hence we reach the
following conclusion:

\begin{Conclusion}\label{conclusion2}

\noindent (1) There is a finite set $\GG^2(M)$ of compact,
orientable, connected, $\partial$-patterned 3-manifolds such that
for each closed, incompressible surface $S\subset M$, all
patterned components of $(G^2_S, G^2_S\cap A^2_S)$ belong to
$\GG^2(M)$;

\noindent (2) Each component of $IB^2_S$ is an $I$-bundle over a
surface $F$ such that $\chi(F)<0$ if $\partial F \ne  \emptyset$.
Moreover $A_S^2$ is incompressible and $\partial$-incompressible.
\end{Conclusion}

\noindent {\bf Step 3.} {\it Comparing the decomposition
$G^2_S\cup _{A^2_S} IB^2_S$ with the GI-decomposition.}

We recall that $M_S=G_S\cup_ {A_S} IB_S^-$ is the
$GI$-decomposition. By the embedded version of the enclosing
property of the JSJ-decomposition and Conclusion
\ref{conclusion2},  $IB^2_S$ is a sub-I-bundle of $IB^-_S$ up to a
proper isotopy of $M\setminus S$. Hence

$$G_S^2=M_S\setminus IB^2_S = G_S\cup_{A_S} (IB_S^-\setminus IB^2_S).$$

Suppose $A_S^{2}$ has $m_S$ components. Let $T^*$ be the once
punctured torus and define the patterned 3-manifold $(P, A)=(T^*\times
I,\partial T^*\times I).$ Let $M_S^{2*}$ be obtained from
$G_S^{2}$ and $m_S$ copies of $(P,A)$ by identifying each frontier
annulus of $G_S^2$ with a frontier annulus of $P$. Then $M_S^{2*}$
is boundary irreducible and is uniquely determined by $G_S^{2}$.
In particular there are finitely many topological types of
$M_S^{2*}$ for all incompressible surfaces $S\subset M$ by
Conclusion \ref{conclusion2} (1).

Let $M_S^{2*}=G^*_S\cup _{A^{*}_S} IB^{*-}$ be the
$GI$-decomposition, which is unique up to isotopy. Hence there are
finitely many topological types of $G_S^{*}$ for all
incompressible surfaces $S\subset M$. It is not difficult to see
that $(G^*_S, A^*_S)=(G_S, A_S)$ for each incompressible surface
$S\subset M$. Hence Theorem 2.3 is proved.\qed

\begin{Definition} Let $M$ be a closed orientable irreducible 3-manifold.
Define:

$\MM=\{M_S, \tilde M_S| \text{ where $S$ runs over all
incompressible surfaces in $M$},$ and $\tilde M_S$ runs over all
double coverings of $M_S \}$.
\end{Definition}

Since each compact 3-manifold has only finitely many double
coverings, the main results in Sections \ref{guts} and
\ref{norm-volume} and their proofs imply the following
proposition:

\begin{Corollary}\label{reduction1}
Let $M$ be a closed, irreducible 3-manifold. Then
there is a finite set $\widetilde \GG(M)$ of
connected compact $\partial$-patterned 3-manifolds such that for
any $X\in \MM$, each component of the patterned guts $(G_X, A_X)$
belongs to $\widetilde \GG(M)$.
\end{Corollary}

\section{Thurston norm and Gromov volume}\label{norm-volume}

\subsection{Finiteness of the Thurston norm}\label{finite norm}

We first give a brief description of the Thurston norm on the
second relative homology group $H_{2}(X, Y;\Z)$ of a compact,
orientable 3-manifold $X$, where $Y \subset \partial X$ is a
subsurface.

Thurston \cite{Th3} introduced a pseudo norm on $H_{2}(X,Y;\Z)$
using the fact that any homology class $z \in H_{2}(X,Y;\Z)$ can be
represented by a properly embedded oriented surface $(F,\partial F)
\hookrightarrow (X, Y)$. Set $\chi_-(F)=\text {max}\{0,- \chi(F)\}$
if $F$ is connected, otherwise let $\chi_-(F)=\sum \chi_-(F_i)$,
where $F_i$ are the components of $F$. Then for an integral class
$z\in H_{2}(X,Y;\Z)$, the \emph{Thurston norm} $\|z\|$ of $z$ is
defined as
$$
\begin{array}{rl}
\|z\| = {\rm inf} \big\{ & \hspace{-3mm}\chi_-(F): F
\mbox{ is an
embedded closed orientable surface } \\
& \hspace{-3mm}\mbox{representing the homology class $z$ in }
H_{2}(X, Y;\Z) \;\big\}.
\end{array}
$$
Thurston then shows that $\|\,\,\|$ extends to a convex
pseudo-norm on $H_{2}(X, Y;\R)$ which is linear on rays through
the origin. The Thurston norm turned out to be very useful in the
study of the topology of 3--dimensional manifolds.

In \cite{Ga1} (see also \cite{Pe}) Gabai shows that to define the
Thurston norm, one can replace ``embedded surfaces'' by ``singular
surfaces'' and still get the same norm.

\begin{Definition} Let $X$ be a compact, orientable 3-manifold and $Y
\subset \partial X$ be a subsurface. For a finite set of elements
$\alpha=\{a_1,..., a_k\}$ of $H_2(X,Y; \Z)$, we define

$$TN(\alpha)=\text{max}\{\|a_i\|, i=1,...,k\}.$$

\noindent Then we define
$$TN(X,Y)=\text{min}\{TN(\alpha) |   \alpha \text{ runs over all
finite sets of elements  of}$$ $$ H_2(X,Y;\Z) \text{which generate}
\, H_2(X,Y;\Q)\}
$$

\noindent to be the Thurston norm of the pair $(X,Y)$.
\end{Definition}

\begin{Lemma}\label{double1}
Let  $p: (\tilde X, \tilde Y)\to (X, Y)$ be a proper non-zero degree map. Then
$TN(\tilde X, \tilde Y)\ge TN(X,Y)$.
\end{Lemma}

\demo Suppose
$\alpha=\{a_1,..., a_k\}\subset H_2(\tilde X,\tilde Y; \Z)$
generates $H_2(\tilde X,\tilde Y; \Q)$. Let $(S_i,\partial
S_i)\subset (\tilde X, \tilde Y)$ be a proper surface which presents
$a_i$ and realizes its Thurston Norm.

Clearly $p(\alpha)=\{p(a_1),..., p(a_k)\}\subset H_2( X, Y; \Z)$.
Since non-zero degree maps induce surjections on rational homology,
$p(\alpha)=\{p(a_1),..., p(a_k)\}$ generates $H_2( X, Y; \Q)$. Now
each $p(S_i)$ is a singular surface representing $p(a_i)$. By
Gabai's result \cite{Ga1}, it follows that $||p(\alpha)||\le ||
\alpha||$, and therefore Lemma \ref{double1} is derived.
 \qed

Recall for each closed incompressible surface $S\subset M$, we have
the $GI$-decomposition $M_S=G_S\cup_ {A_S} IB^-_S$.

\begin{Lemma}\label{double2}
There is a double cover $\tilde M_S=\tilde G_S\cup_ {\tilde A_S}
\widetilde {IB^-_S}$ of $M_S=G_S\cup_ {A_S} IB^-_S$ such that
each component of  $\widetilde {IB^-_S}$ is a product of an
orientable surface with the interval.

\end{Lemma}

\demo An elementary fact is that each compact non-orientable surface
$F$ is doubly covered by an orientable surface $\tilde F$ such that
the restriction on each component of $\partial \tilde F$ is an
homeomorphism. It follows that each twisted $I$-bundle $B$ over
compact non-orientable surface $F$ is doubly covered by a product
$I$-bundle $\tilde B=\tilde F\times I$ such that the restriction on
each component of $\partial \tilde F\times I$ is a homeomorphism.

For each twisted $I$-bundle component $B$ of  $(IB^-_S, A_S)$, pick
a double covering given in the first paragraph, and for each
remaining component of  $(IB^-_S, A_S)$ and each component of $(G_S,
A_S)$, pick two identical patterned copies of it. Obviously we can
glue them together to get a double cover $p:\tilde M_S\to M_S$. Let
$\tilde G_S$, ${\tilde A_S}$ and $\widetilde {IB^-_S}$ be the
pre-images of $G_S$, ${A_S}$ and $IB^-_S$, then one can verify from
the definitions that $\tilde M_S=\tilde G_S\cup_ {\tilde A_S}
\widetilde {IB^-_S}$ is the $GI$-decomposition of $\tilde M_S$ and
verifies the desired property. \qed

Now we are going to prove the main result of this section.

\begin{Theorem}[Finiteness  of the Thurston norm]\label{norm}
Let $M$ be an irreducible, closed, orientable 3-manifold. Then
$TN(M_S,
\partial M_S)$ takes
at most finitely many values, when $S$ runs over all closed,
incompressible surfaces embedded in $M$.
\end{Theorem}

\demo By Lemmas \ref{double1}, we need only to prove Theorem
\ref{norm} for double coverings $(\tilde M_S, \tilde \partial M_S)$
provided by Lemma \ref{double2}, for all incompressible surfaces $S
\subset M$. For simplicity we still use $M_S=G_S\cup_ {A_S} IB^-_S$
to denote $\tilde M_S=\tilde G_S\cup_ {\tilde A_S} \widetilde
{IB^-_S}$. Then by Corollary \ref{reduction1} there are only finitely
many topological types of patterned guts $(G_S, A_S)$ for all
incompressible surfaces $S \subset M$. Hence the number of
components of $A_S$ is uniformly bounded. Again by Lemma
\ref{double2},  each component of ${IB^-_S}$ is a product of an
orientable surface with the interval.

We first modify the decomposition so that the gluing annuli between
the two parts become separating. For each component $N(F)$ of
$IB^-_S$ we choose a curve in the interior of the base surface $F$,
which co-bounds a planar subsurface $Q$ together with all the
boundary components of $F$. Since the number of boundary components
of $F$ is bounded by the number of components of $A_S$, $|\chi(Q)|$
is uniformly bounded above, for all incompressible surfaces $S
\subset M$. Then we consider the new decomposition $M = G'_S
\cup_{A'_S} {IB'}^-_S$, where $G'_S$ is obtained by gluing to $G_S$
the handlebodies $N(Q)$ along the components of $A_S$, and
${IB'}^-_S$ is the sub-$I$-bundle of $IB^-_S$ corresponding to the
subsurfaces $F - \text{int}(Q)$. The gluing annuli $A'_S$ are the
separating annuli of $N(\partial Q)- A_S$, using our previous
convention that $N(Q)$ and $N(\partial Q)$ are the $I$-bundle
restricted to $Q$ and $\partial Q$ respectively.

For a given patterned guts $(G_S, A_S)$, there are only finitely
many positive integer solutions $\{m_1,...., m_k\}$ such that
$m_1+...+m_k=m$ where $m$ is the number of components of $A_S$,
and for any such solution $\{m_1,..., m_k\}$, there are only
finitely many ways to distribute $m$ elements into $k$ groups of
cardinality $m_1,...,m_k$ respectively. Hence by the construction
and Theorem \ref{guts finiteness},  there are only finitely many
topological types of $\partial$-patterned 3-manifolds $(G'_S,
A'_S)$ for all incompressible surfaces $S$ in $M$. Then the
finiteness for the values of $TN(M_S, \partial M_S)$ is a direct
consequence of the following lemma:

\begin{Lemma}\label{basis norm} Let $S\subset M$ be a closed incompressible
surface, then:

$TN(M_S,\partial M_S)\leq TN(G'_S, \partial G'_{S} \setminus
\text{int}A'_S)$.
\end{Lemma}

\demo We consider the following natural homomorphisms induced by
the inclusion maps:

\noindent $\phi: H_2(G'_S,\partial G'_{S} \setminus
\text{int}A'_{S};\Z) \to H_2(M_S,\partial M_S;\Z)$;

\noindent$\psi: H_2({IB'}^-_S,\partial {IB'}^-_S \setminus
\text{int}A'_{S};\Z) \to H_2(M_S,\partial M_S;\Z)$.

By applying the relative Mayer-Vietoris sequence (see \cite[page
52]{Do}) to the pairs $(G'_S,\partial G'_{S} \setminus
\text{int}A'_{S})$ and $({IB'}^-_S,\partial {IB'}^-_S \setminus
\text{int}A'_{S})$, one gets the exact sequence:

\medskip

$ \dots \to H_2(A'_S, \partial A'_S; \Z) \to H_2(G'_S,\partial
G'_{S} \setminus \text{int}A'_{S};\Z) \oplus H_2({IB'}^-_S,\partial
{IB'}^-_S \setminus \text{int}A'_{S};\Z) $

\hskip 0.5 true cm \noindent $\to H_2(M_S,\partial M_S;\Z) \to
H_1(A'_S,
\partial A'_S; \Z)$

\hskip 0.5 true cm \noindent $\to H_1(G'_S,\partial G'_{S} \setminus
\text{int}A'_{S};\Z) \oplus H_1({IB'}^-_S,\partial {IB'}^-_S
\setminus \text{int}A'_{S};\Z)\to\dots $

We first show the injectivity of the homomorphism
$$H_1(A'_S,
\partial A'_S; \Z)\to H_1(G'_S,\partial G'_{S} \setminus \text{int}A'_{S};\Z) \oplus
H_1({IB'}^-_S,\partial {IB'}^-_S \setminus \text{int}A'_{S};\Z).$$
To do this we need only to show the
injectivity of each homomorphism
$$H_1(A, \partial A; \Z)\to
H_1({IB'}^-_A,\partial {IB'}^-_A \setminus \text{int}A;\Z),$$ where
$A$ is a component of $A'_S$, and ${IB'}^-_A$ is the component of
${IB'}^-_S$ containing $A$.

Note $H_1(A,
\partial A; \Z)=\Z$ is generated by any arc in $A$ connecting the two
components of $\partial A$, and ${IB'}^-_A=F\times [0,1]$, where $F$
is an orientable surface with $\partial F\times [0,1]=A$. Let $F^*$
be a proper oriented surface  of $\partial F\times [0,1]$. Then it
is a direct geometric observation that the number of  times that
$\partial F^*$ crosses $A$ from $\partial F\times \{0\}$ to
$\partial F\times \{1\}$ and from $\partial F\times \{1\}$ to
$\partial F\times \{0\}$ must be the same. This shows the required
injectivity.

Then by the exact sequence we have $\partial_*:  H_2(M_S,\partial
M_S;\Z) \to H_1(A'_S,
\partial A'_S; \Z)$ is null, and thus we get an epimorphism:

\noindent $\phi + \psi :H_2(G'_S,\partial G'_{S} \setminus
\text{int}A'_{S};\Z) \oplus H_2({IB'}^-_S,\partial {IB'}^-_S
\setminus \text{int}A'_{S};\Z) \to H_2(M_S,\partial M_S;\Z)$.

It is clear that $H_2({IB'}^-_S,\partial {IB'}^-_S \setminus
\text{int}A'_{S};\Z)$ has a basis $\gamma = \{c_1,\dots , c_m\}$
which is formed by a set of vertical annuli, whose Thurston norm
vanishes. Hence for any generating set $\beta = \{b_1, \dots ,
b_n\}$ of $H_2(G'_S,\partial G'_{S} \setminus
\text{int}A'_{S};\Z)$, $\alpha = \{\phi(b_1), \dots , \phi(b_n),
\psi(c_1), \dots , \psi(c_m)\}$ is a generating set of
$H_2(M_S,\partial M_S;\Z)$. It follows that $TN(M_S, \partial M_S)
\leq TN(\alpha) \leq TN(\beta)$, since by the definition of
Thurston norm $\|\phi(b_i)\| \leq \|b_i\|$ and $0 \leq
\|\phi(c_j)\| \leq \|c_j\| = 0$, for $i=1,...,n$ and $j=1,...,m$.
Therefore $TN(M_S,
\partial M_S) \leq TN(G'_S, \partial G'_{S} \setminus \text{int}A'_S)$.\qed

\subsection{Finiteness of absolute Gromov volume}\label{finite volume}

This section will not be used in the rest of the paper, but it
provides a finiteness result for the absolute Gromov volumes of
the compact manifolds $M_S$, analogous to the one for their
Thurston norms.

First we recall the basic definitions about Gromov's simplicial
volume (see [Gr]).

\begin{Definition} Let $X$ be a compact orientable 3-manifold with boundary.
Define the \emph{relative Gromov volume} $\vert X,\partial X
\vert$ by:

\[
{\vert X,\partial X \vert}:= \inf \left\{ \sum_{i=1}^n
\vert\lambda_i\vert \left \vert
    \begin{array}{l}
        \sum\limits_{i=1}^n \lambda_i\sigma_i \text{ is a cycle representing a
fundamental} \\
        \text{class in } H_3(X,\partial X;\R),\text{ where} \,
\sigma_i:\Delta^3\to X \\
        \text{is a singular simplex and } \lambda_i\in\R,\ i=1,\ldots,n.
    \end{array}
\right.\right\}
\]
\end{Definition}

A fundamental class in $H_{3}(X,\partial X;\R)$ is the image of
any of the $2^k$ fundamental classes in $H_{3}(X,\partial X;\Z)$
under the coefficient homomorphism, where $k$ is the number of
connected components of $X$.

For a manifold with non-empty boundary, there is another way of
defining a simplicial volume, that we call the \emph{absolute
Gromov volume}.

\begin{Definition} Let $X$ be a compact orientable 3-manifold with boundary and let
$D(X)$ be the double of $X$, obtained by identifying two copies of
$X$ along their boundary via the identity map. The absolute Gromov
volume of $X$, denoted as $\vert X \vert$, is defined to be half
of the Gromov volume of the closed manifold $D(X)$.
\end{Definition}

By the definitions of these two volumes, one has: $\vert X \vert
\leq \vert X,\partial X \vert$. Moreover by [So5] and [Ku] they
are equal if and only if $\partial X = \emptyset$ or
$\chi(\partial X) \geq 0$.

For example, let $(X,A)$ be a patterned 3-manifold and let
$D_A(X)$ be the compact 3-manifold obtained by doubling $X$ along
the portion $\partial X\setminus A$ of its boundary. Since
$\partial D_A(X)$ is a collection of tori, one has $\vert D_A(X)
\vert = \vert D_A(X),\partial D_A(X) \vert$.

The following finiteness result holds for the absolute Gromov
volume, while it is false for the relative Gromov volume.

\begin{Proposition}[Finiteness of Gromov volume]\label{volume}
Let $M$ be an irreducible, compact, orientable 3-manifold. Then
the absolute Gromov volume $\vert M\setminus S \vert$ takes only
finitely many values for all incompressible surfaces $S\subset M$.
\end{Proposition}

\demo For an incompressible surface $S\subset M$, we consider the
$GI$-decomposition $M_S=G_S\cup_ {A_S} IB^-_S$. By Conclusion
\ref{conclusion2} (2) in the proof of Theorem \ref{guts finiteness},
$\partial D(A_S)$ is a collection of incompressible tori in
$D(M_S)$. Since the relative Gromov volume is additive under gluing
along incompressible tori (see \cite{So5}), we have:

\noindent $|D(M_S)|=|D_A(G_S), \partial D_A(G_S)|+|D_A(IB^-_S),
\partial D_A(IB^-_S)|$.

Moreover the relative Gromov volume of $D_A(IB^-_S)$ vanishes,
because $D_A(IB^-_S)$ is homeomorphic to an $S^1$-bundle (see
\cite{Gr}, \cite[Chap. 6]{Th1},) and the relative Gromov volume of
$D_A(G_S)$ equals its absolute Gromov volume. Therefore we get
$\vert M_S \vert = \frac{1}{2} \vert D_A(G_S) \vert$.

Now Theorem \ref{guts finiteness} shows that there are only
finitely many possible topological type for $D_A(G_S)$ when $S$
runs over all incompressible surfaces $S$ in $M$, and hence
Proposition \ref{volume} follows.\qed

\begin{Example} We give an example of a closed orientable 3-manifold $M$
such that the relative Gromov volume $\vert M_S,\partial M_S \vert$
is unbounded when $S$ runs over all incompressible surfaces $S$ in
$M$. Let $X$ be a knot exterior in $\S^3$ which contains
incompressible Seifert surfaces of arbitrarily high genus (such
examples exist and can even be hyperbolic see \cite{Gu}). Then the
closed manifold $M = D(X)$ contains non-separating incompressible
closed surfaces $S_n$, with $\chi_-(S_n)$ tending to infinity with
$n$, formed by doubling the Seifert surfaces. When $M$ is split open
along such  a surface, by definition of the relative Gromov volume,
one has: $\vert M_{S_n},
\partial M_{S_n} \vert \geq 2 \chi_- (\partial M_{S_n}) = 4\chi_-(S_n)$. Hence  $\vert M_{S_n}, \partial M_{S_n} \vert$ tends
to infinity with $n$.
\end{Example}

\section{Local Domination}\label{JSJ-pieces}

In this section we prove the finiteness of the JSJ-pieces for
 manifolds which are dominated by a given compact,
orientable 3-manifold. We recall the statement that we are going
to prove:

\begin{Theorem}[Finiteness of JSJ-pieces]\label{local domination}
Let $M$ be a closed, orientable, 3-manifold. Then there is a
finite set $\HH\SS(M)$ of complete hyperbolic 3-manifolds with
finite volume and of Seifert fibered 3-manifolds, such that the
JSJ-pieces of any closed, orientable, irreducible
3-manifold $N$ dominated by $M$ belong to $\HH\SS(M)$, provided
that  $N$ is not supporting the geometries of $\S^3$,
$\widetilde{PSL_2 (\R)}, Nil.$
\end{Theorem}

By \cite[Prop. 3.3]{BW}, we can find an irreducible (even
hyperbolic) closed, orientable 3-manifold which 1-dominates $M$.
Hence in the remainder of the proof, we may assume that $M$ is
irreducible.

Let $M$ be a closed orientable irreducible 3-manifold. By Haken's
finiteness theorem, there is a maximum number $h(M)$ of pairwise
disjoint, non-parallel, closed, connected, incompressible surfaces
embedded in $M$. The following elementary fact (see \cite{W1} for
example) will be used in this section and the next ones.

\begin{Lemma}\label{haken} Let $M$ and $N$ be two closed, irreducible
and orientable 3-manifolds. If $M$ dominates $N$, then $h(M)\geq
h(N)$.
\end{Lemma}

Let  $\Gamma(N)$ be the dual graph associated with the
JSJ-decomposition of $N$. This graph has one vertex for each
Seifert piece or piece with a hyperbolic metric of finite volume
and one edge for each incompressible torus boundary component of
either type of piece. If $M$ dominates $N$, then $h(M)$ gives an
upper bound for the number of edges of $\Gamma(N)$, by Lemma
\ref{haken}. Hence the number of JSJ-pieces of $N$, which is the
number of vertices of $\Gamma(N)$, is bounded above by $h(M)+1$.
Therefore to prove Theorem \ref{local domination}, we need only
show that the JSJ-pieces of all 3-manifolds $N$ dominated by a
closed, orientable 3-manifold $M$ admit  only finitely many
topological types.

Recall the definition: $\MM=\{M_S, \tilde M_S| \text{ where $S$ runs over all
incompressible surfaces in $M$},$ and $\tilde M_S$ runs over all
double coverings of $M_S \}$.

By the proof of Theorem \ref{norm}, we have:

\begin{Corollary}\label{reduction3}
$Sup \{TN(X,\partial X)| X\in\MM \} \leq L(M)$ for
some constant $L(M)>0$ depending only on $M$.\qed
\end{Corollary}

\begin{Proposition}\label{reduction2}
For a given integer $L>0$, there is a finite set $\SS(L)$ of
compact Seifert 3-manifolds such that if a Seifert manifold $N$
with non-empty boundary and orientable base is dominated  by a
compact orientable 3-manifold $P$ with $TN(P,\partial P) \leq L$,
then $N$ belongs to $\SS(L)$.
\end{Proposition}

\demo Each homology class $y$ of $H_2(N,\partial N; \Z)$ can be
represented by an orientable incompressible and
$\partial$-incompressible surface. Since $N$ is an irreducible
Seifert manifold, each incompressible surface is properly isotopic
to either a vertical torus or annulus (foliated by Seifert
circles), or a horizontal surface (transverse to all Seifert
circles) (cf. [Ja, Chap. VI]). Since $\partial N\ne \emptyset$,
$N$ always admits horizontal surfaces.

Let $\OO$ be the orbifold base of $N$ and $h$ be a regular fiber
of $N$. Suppose also that $\OO$, $h$ and $N$ are compatibly
oriented. Let $F$ be a horizontal surface of $N$ and $p: F\to \OO$
the branched covering, induced by the restriction to $F$ of the
projection of $N$ onto its base. Since $\OO$ is oriented, so is
$F$. Note that the Euler characteristic $\chi(\OO)$ is computed
for an orbifold, so that each exceptional fiber of multiplicity
$n$ gives a term ${\frac 1 n} -1$. Then we have $\chi(F)=|d|\times
\chi(\OO)<0$, where $d  =\text{deg}(p) \not = 0$ is equal to the
algebraic intersection number $[F]\cdot [h]$ of $F$ and $h$. Up to
reversing the orientation of $F$, one can always assume that $d =
[F]\cdot [h] >0$. Note that the geometric intersection number
$|F\cap h|$ of $F$ and $h$ (i.e. the minimal number of
intersection points between $F$ and $h$ up to ambient isotopy) is
precisely the absolute value of the algebraic intersection number
$[F]\cdot [h]$.

Suppose further that $F$ has minimal genus among all horizontal
surfaces in $N$. If $\chi(F) \geq 0$, $F$ is a disc or an annulus,
and thus $N$ can be homeomorphic only to a solid torus, an
$S^1$-bundle over the annulus or a twisted $S^1$-bundle over a
M\"obius band. Assuming that these three Seifert manifolds belong
to $\SS(L)$, we can suppose  furthermore that $\chi(F) < 0$.

Let $\|\,\,\|_{N}$ (resp. $\|\,\,\|_{P}$) be the Thurston norm on
$H_2(N,\partial N; \Z)$ (resp. $H_2(P,\partial P; \Z)$). Note that
$H_2(N,\partial N; \Z)$ is torsion free and therefore it is
precisely the integer lattice of $H_2(N,\partial N; \R)$. Let $V =
\{ y \in H_2(N,\partial N; \Z); \| y \|_{N} = 0 \}$. By the
discussion above, $V$ is the sublattice of $H_2(N,\partial N; \Z)$
generated by the vertical tori and annuli.

\begin{Lemma}\label{direct sum}
$H_2(N,\partial N; \Z) = \langle [F]\rangle  \oplus V$.
\end{Lemma}

\demo Pick any homology class $y\in H_2(N, \partial N; \Z)$. If
$\| y \|_{N} = 0$, then $y\in V$. Suppose $\| y \|_{N} \ne 0$. Let
$S$ be an orientable, incompressible and $\partial$-incompressible
surface representing $y = [S]$ with $-\chi(S)= \| y \|_{N}$. Since
$\chi(S)= - \| y \|_{N} <0$, after a proper isotopy we may assume
that $S$ is horizontal and that $[S].[h] = \vert S\cap h \vert
>0$ (otherwise we replace $y$ by $-y$). Let $\ell \geq 1$ be the
integral part of $[S].[h]/[F].[h]$. Then $(\ell +1) ([F].[h])>
[S].[h] \geq \ell ([F].[h])$, that is $|F\cap h|= [F]\cdot [h] >
[S-\ell F]\cdot [h] \geq 0.$

If the homology class $[S-\ell F]$ does not belong to $V$, then it
can be represented by a horizontal surface $S'$ such that :

\noindent$\| [S-\ell F] \|_{N} = -\chi(S') = -([S-\ell F]\cdot
[h])\chi(\OO) > -[F]\cdot [h]\chi(\OO) = \| [F] \|_{N}$.

This would contradict the minimality of the genus of $F$ among all
horizontal surfaces in $N$. Therefore $[S-\ell F]\in V$ and
$y=[\ell F]+[S-\ell F]$.\qed

By hypothesis, there is a compact, orientable 3-manifold $P$ with
$TN(P,\partial P) \leq L$ and a non-zero degree map $f: P \to N$.
Let $\alpha =\{z_1,...,z_m \}$ be a basis of $H_2(P,\partial P;
\Z)$ realizing  $TN(P)$ : $\text{max}\{\|z_i\|_{P}; i=1,...,m\}
\leq L$. For $i=1,...,m$, let $S_i$ be a properly embedded surface
in $P$ representing $z_i$ with $-\chi (S_i) = \|z_i\|_{P}$.

For $i=1,...,m$, we set $y_i=[f(S_i)]=\ell_{i} [F] + v_{i} \in
H_2(N,\partial N; \Z)$, where $v_{i}\in V$. By the triangle
inequality and the fact that $\|v_i\|_{N} = 0$, we get:

$$\vert \ell_{i} \vert \|[F]\|_{N} = \|\ell_{i} [F]\|_{N} = \|y_i -
v_{i}\|_{N} \leq \|y_i\|_{N} + \|v_i\|_{N} = \|y_i\|_{N}$$

By [Ga1] (see also [Pe]) the Thurston norm $\|y_i\|_{N}$ can be
calculated using singular surfaces, therefore $\|y_i\|_{N} \leq
-\chi (S_i) = \|z_i\|_{P} \leq L$. Combining the two inequalities,
we have $\ell_{i} \|[F]\|_{N} \leq L$ for $i=1,...,m$.

Since $f: P \to N$ has non-zero degree, $f_*(H_2(P,\partial P;
\Z))$ has finite index in $H_2(N,\partial N; \Z)$ and thus it
cannot lie in $V$. Therefore there is some index $i \in \{ 1,...,m
\}$ with $\vert \ell_i \vert \geq 1$. It follows that $\|[F]\|_{N}
\leq L$, hence the horizontal surface $F$ can have only finitely
many topological types, up to homeomorphism.

Cutting the Seifert manifold $N$ along the horizontal surface $F$,
we obtain a product $F \times I$, since the base $\OO$ and the
surface $F$ are orientable. Therefore $N$ can be presented as a
surface bundle over $S^1$ with fiber $F$ and orientation preserving
monodromy $g:F\to F$. Since $N$ is Seifert fibered, $g$ must be a
periodic map \cite[Chap. VI]{Ja}. However, up to conjugacy, a given
compact surface  admits only finitely many periodic homeomorphisms.
Since any two conjugate monodromy maps define homeomorphic
3-manifolds, there are only finitely many possible homeomorphism
types of Seifert manifolds $N$ for a given compact surface $F$.
Since $F$ has only finitely many topological types, the proof of
Proposition \ref{reduction2} is complete.\qed

\medskip

\noindent{\bf Proof of Theorem \ref{local domination}.}

We assume first that $N$ is not a Seifert manifold. By Soma's
results (\cite{So1}, \cite{So2}), we know that Theorem \ref{local
domination} holds for hyperbolic JSJ-pieces. Since $N$ is not
Seifert fibered, we may assume that $N$ has a non-empty JSJ-family
of tori $\TT$ and we have only  to consider the Seifert fibered
JSJ-pieces.

Let $f: M\to N$ be a map of non-zero degree. After a homotopy of
$f$ we may assume that $f^{-1}(\TT)$ is a non-empty collection of
disjoint non-parallel closed incompressible surfaces in $M$. Let
$\MM_f$ be the union of all components of $M\setminus f^{-1}(\TT)$
and of all their double coverings. By definition we have $\MM_f
\subset \MM$.

Let $N_i \subset N$ be a Seifert fibered JSJ-piece. Then $N_i$ is
dominated by at least one component $M_i$ of $M\setminus
f^{-1}(\TT)$. Then the finiteness of such JSJ-pieces $N_i$ with an
orientable base follows immediately from
Corollary \ref{reduction3}  and Proposition \ref{reduction2}.

If $N_i$ has a non-orientable base orbifold, let $\tilde N_i$ be
the unique double cover of $N_i$ which is Seifert fibered with an
orientable base. Then a standard argument shows that a double cover
$\tilde M_i$ of $M_i$ dominates  $\tilde N_i$. Thus we get the
finiteness of such 3-manifolds $\tilde N_i$ as above from Corollary
 \ref{reduction3} and Proposition \ref{reduction2}. Since any
involution on such Seifert manifolds $N_i$ is conjugate to a fiber
preserving one by \cite{MS}, there are only finitely many conjugacy
classes of involutions on each $\tilde N_i$. This implies the
finiteness of the Seifert JSJ-pieces $N_i$.

The finiteness of Seifert manifolds $N$ supporting a product
geometry $\H^2 \times \R$ follows also from Corolary
\ref{reduction1} and Proposition \ref{reduction2} as above, and thus
Theorem \ref{local domination} is proved (see also \cite{WZ}).\qed

Using a standard doubling construction, Theorem \ref{local
domination} can be extended to the following case where the
3-manifold targets have toric boundary.

\begin{Corollary}\label{local domination with boundary} Let $M$ be a compact, orientable, 3-manifold. Then there is a
finite set $\HH\SS(M)$ of complete hyperbolic 3-manifolds with
finite volume and of Seifert fibered 3-manifolds, such that the
JSJ-pieces of any compact, orientable, irreducible, 3-manifold $N$
with non-empty toric boundary, dominated by $M$ belong to
$\HH\SS(M)$.
\end{Corollary}

\demo If $N$ has non-empty boundary, so does $M$. Let $D(N)$ be the
double of $N$, obtained by gluing two copies of $N$ along their
boundaries via the identity map . Then the double $D(M)$ of $M$
dominates $D(N)$. Since the boundary of $N$ is a collection of
tori, the JSJ-pieces of $D(N)$ are either exactly those of $N$ and
consist of two copies of hyperbolic and Seifert pieces in the
JSJ-decomposition of $D(N)$, or there are new Seifert fibered
pieces obtained by doubling some Seifert fibered pieces of $N$
along some of their boundary tori. In any case, the finiteness of
the JSJ-pieces of $D(N)$ implies the finiteness of the JSJ-pieces
of $N$. \qed

\begin{Remark}\label{compact} A similar double construction argument shows that finiteness of closed
targets implies finiteness of the targets in the setting of compact
orientable 3-manifolds. First finiteness of irreducible and $\partial$-irreducible compact targets implies finiteness of
 compact targets. Since the double $D(M)$ of an  irreducible and $\partial$-irreducible compact 3-manifold $M$ is
 Haken, there are, up to conjugacy, only finitely many involutions with 2-dimensional fixed point set on $D(M)$
 by \cite{To} and the proof of the geometrization conjecture for Haken manifolds. Therefore only finitely many irreducible
 and $\partial$-irreducible compact 3-manifolds have homeomorphic doubles.
\end{Remark}

\section{Integral homology spheres}\label{spheres}

The main result of this section gives a positive answer to
Question 3 when the targets are integral homology spheres. It
implies a positive answer to Question 2 when the targets are
integral homology spheres and to Question 1 when the
domain is an integral homology sphere.

\noindent {\bf Theorem 1.2} \emph{Any closed orientable $3$-manifold
dominates at most finitely many integral homology spheres.}

Let us fix $M$ as a closed orientable $3$-manifold. As in the
previous section, we may assume for the remainder of the proof
that $M$ is irreducible.

First we reduce the proof to the case where the target homology
sphere $N$ is irreducible. As in the previous section, the
preimage of a collection of separating essential spheres
associated with the prime decomposition of $N$ can be assumed to
be incompressible, disjoint and non-parallel surfaces in $M$.
Hence  there are at most $h(M) + 1$ prime factors. Moreover by
pinching all the prime factors except one to a point, it follows
that each prime factor is dominated by $M$. Hence we have only to
show the finiteness of the set $\DD(M)$ of homeomorphism classes
of irreducible, integral homology
$3$-spheres $N$ which are dominated by $M$.

A {\it slope} on a torus $T$ is an isotopy class of essential simple
closed curves. The set of slopes on $T$ corresponds bijectively with
$\pm$-classes of primitive elements of $H_1(T; \mathcal{Z})$.

Given a slope $\alpha$ on a torus boundary component $T$ of a
$3$-manifold $Y$, the {\it $\alpha$-Dehn filling} of $Y$ with slope
$\alpha$ is the $3$-manifold $Y(\alpha) := (S^1 \times D^2) \cup_f
Y$ where $f$ is any homeomorphism $\partial(S^1 \times D^2) \to T$
such that $f(\{*\} \times \partial D^2)$ represents $\alpha$. It is
well-known that $Y(\alpha)$ is independent of the choice of $f$.

First let us recall some definition and primary facts about Seifert
manifolds.

 Let  $F_{g,n}$ be an
oriented  $n$-punctured surface of genus $g$ with boundary
components $c_1,...,c_n$ with $n\ge 0$. Then $N'=F_{g,n}\times S^1$
is oriented if $S^1$ is oriented. Let  $h_i$ be the oriented $S^1$
fiber on the torus $c_i\times h_i$
 (call such pairs $\{(c_i, h_i)\}$ a section-fiber coordinate
system). Let $0\leq s\leq n$, we attach $s$ solid tori $V_i$ to the
boundary tori of $N'$ such that the meridian of $V_i$ is identified
with the slope $r_i=c_i^{a_i}h_i^{b_i}$ where $a_i>0, (a_i, b_i)=1$
for $i=1,...,s$. We denote $N(
g,n-s;\frac{b_1}{a_1},\cdots,\frac{b_s}{a_s})$ the resulting
manifold which has the Seifert fibered structure extended from the
circle bundle structure of $N'$. Each orientable Seifert fibered
space with orientable base $F_g$ with $n-s$ boundary components and
$s$ exceptional fibers is obtained in such a way.

\begin{Lemma}\label{Seifert-Basic}
Suppose $N$ is a Seifert manifold given as above.

(1) Suppose $N$ is an integer homology 3-sphere. Then $N$ is closed
and $g=0$, furthermore $({\prod_{i=1}^na_i})(\sum_1^n
\frac{b_i}{a_i})=1$.

(2) Suppose $n>s$ and $N(\mu_{s+1}, ..., \mu_n)$ is an integer
homology 3-sphere, where $\mu_j=(a_i, b_j)$, $j= s+1, ..., n$, then $g=0$ and moreover

(i) if each $a_j>0$ for $j\in \{s+1, ..., n\}$, then the Seifert
fibration of $N$ extends over $N(\mu_{s+1}, ..., \mu_n)$ and
$(\prod_{i=1}^na_i)(\sum_1^n \frac{b_i}{a_i})=1$;

(ii) if some $a_j=0$ for some $j\in \{s+1, ..., n\}$, then
$b_j\prod_{i=1, i\ne j}^{ n} a_i=1$.
\end{Lemma}

\demo The proof of the lemma is an application of linear algebra.
(1) is well known, see \cite[3.1]{HWZ1} for example. (2) (i) mainly
follows from (1).

For (2) (ii), if $a_j=0$ for some $j\in \{s+1, ..., n\}$, then $b_j$
must be 1. For some $i\ne j$, $a_i=0$ implies that $N(\mu_{s+1},
..., \mu_n)$ has positive first Betti number, and $a_i>1$ implies
that $H_1(N(\mu_{s+1}, ..., \mu_n),Z)$ contain a torsion element of
order $a_i$. So $a_i=1$ for $i\ne j$.\qed

Below we denote $N(0,n-s;\frac{b_1}{a_1},\cdots,\frac{b_s}{a_s})$ as
$N(n-s;\frac{b_1}{a_1},\cdots,\frac{b_s}{a_s})$, and $N(
0,0;\frac{b_1}{a_1},\cdots,\frac{b_s}{a_n})$  as
$N(\frac{b_1}{a_1},\cdots,\frac{b_s}{a_s})$.

\begin{Lemma}\label{seifert} Only finitely many Seifert
fibered integral homology 3-spheres belong to $\DD(M)$.
\end{Lemma}

\demo A Seifert fibered integral homology 3-sphere must
support the geometry of either $S^3$ or $\widetilde
{PSL_{2}(\R)}$.

For Seifert manifolds supporting the geometry $S^3$, there are only
two integral homology 3-spheres: the 3-sphere $S^3$ and the
Poincar\'e dodecahedral space.

Now suppose that $N$  supports the geometry of $\widetilde
{PSL_{2}(\R)}$. Since $N$ is an integral homology sphere, By
\cite{Sc} , $N= (\frac{b_1}{a_1}, \dots , \frac{b_n}{a_n})$ which
satisfy:\begin{itemize}

\item The rational Euler number $e=
-\sum_{i=1}^n\frac {b_i}{a_i}$  is non-zero.

\item The Euler characteristic of the orbifold base B is $\chi(B)
=2 -\sum_{i=1}^n\left(1-\frac 1{a_i}\right)<0$.

\item $\left|e\prod_{i=1}^n a_i\right|=1$ by Lemma
\ref{Seifert-Basic}.
\end{itemize}

Thus $e= \frac {\pm 1}{\prod_{i=1}^na_i}$ and $b_{i}(\prod_{j\not
=i}a_j) \equiv \pm 1 \, \text{modulo} \, a_i$ for $i = 1, \dots ,
n$. Moreover the integers $a_i, i = 1, \dots , n$ are pairwise
relatively prime. Therefore the unordered set $\{a_1, \dots , a_n\}$
of integers determines the Seifert fibered homology  sphere $N$, up
to orientation.  So we need only to show that if
$N$ is dominated by $M$, then $n$ and the integers $a_i, i = 1,
\dots , n$, take only finitely many values. In fact to do so, it is
sufficient to get a uniform upper bound on $\prod_{i=1}^n a_i$, depending
only on $M$.

We use the Seifert volume $SV$ introduced by Brooks and Goldman
\cite{BG}. It has the following interesting properties:

\noindent (1) $SV(M) \geq d SV(N)$ if $f: M\to N$ is a  map of
degree $d \not = 0$, for orientable 3-manifolds $M$ and $N$.

\noindent (2) $SV(N)=\left|\frac {\chi(B)^2}{e(N) }\right|$ if $N$
is a $\widetilde {PSL_2(\R)}$-manifold with base orbifold $B$.

It is easy to see that the maximum of the Euler characteristic of
the base $B$ of $N$ is obtained for the sphere with three cone points with orders
$\{2, 3, 7\}$. Hence:

\noindent (3) $\chi(B) \leq -\frac{1}{42}$.

Then by (1), (2) and (3) we have : $SV(M)\ge d SV(N) = d
\left|\frac {\chi(B)^2}{e(N) }\right|\ge \left|\frac{1}{42^2}
\prod_{i=1}^n a_i\right|$. Therefore  $\prod_{i=1}^n a_i \leq
42^{2} SV(M)$ and the proof of Lemma \ref{seifert} is
complete.\qed

\bigskip

The dual graph $\Gamma(N)$ to the JSJ-decomposition of an
irreducible homology sphere $N$ is a tree. By Lemma \ref{haken},
the number of edges of $\Gamma(N)$ is $\leq h(M)$, the Haken
number of $M$.

By the local domination theorem and Lemma \ref{seifert}, the
geometric JSJ-pieces of the closed orientable $3$-manifolds in
$\DD(M)$ belong to a finite set $\HH\SS(M)$ of compact
$3$-manifolds with interiors admitting complete hyperbolic metrics
with finite volume and of Seifert $3$-manifolds.

For a given graph $\Gamma$, let $\DD(M, \Gamma) \subset \DD(M)$ be
the set of homeomorphism classes of  closed
orientable integer homology $3$-spheres $N$ such
that:\begin{enumerate}

\item $N$ is dominated by $M$.

\item The JSJ-graph $\Gamma(N)$ is abstractly isomorphic to
$\Gamma$.

\item Each vertex manifold has a fixed topological type. Each torus
boundary component of the vertex manifold is assigned to an edge on
the vertex.

\end{enumerate}

The local domination theorem and Lemma \ref{haken} reduce  Theorem
\ref{homology spheres} to the following:

\begin{Proposition}\label{fixed graph}
The set $\DD(M, \Gamma)$ is finite.
\end{Proposition}

Before starting the proof of this proposition, we need to introduce some
notions, definitions and constructions which will be useful.

For each integral solid torus $V$, the kernel of the induced
homomorphism $H_1(\partial V; \Z) \to  H_1(V; \Z)$ is infinite
cyclic, generated by an essential simple loop  which bounds a
properly embedded surface $F_V$ in $V$. Such a surface $F_V$ is
called a Seifert surface for the integral homology solid torus $V$.
Then the slope $\lambda_V \in H_1(\partial V; \Z)$ of $\partial F_V$
does not depend of the Seifert surface and is uniquely determined by
the topological type of $V$. We call $\lambda_V$ the {\it
longitudinal slope} of $V$ on $\partial V$.

\noindent {\bf The pinch construction.} Let $Y$ be a compact,
orientable 3-manifold with $T \subset \partial Y$ a torus boundary
component. Let $Z$ be an integral homology solid torus, $\phi:
\partial Z \to T \subset
\partial Y$ a gluing map  and $Y'=Z\cup_\phi Y$. By pinching a
Seifert surface $F_Z$ onto a disk $D^2$, one can define a proper
degree-one map $p_Z : Z \to S^1 \times D^2$ such that
$p^{-1}_{Z}(\{x\}\times\partial D^2) = \lambda_Z$ for some point $x
\in S^1$. Then one gets a degree-one map $f_Z : Y' \to Y({\mu})$,
which is the identity on $Y$ and where $Y({\mu})$ is obtained by
Dehn filling the component $T$ of $\partial Y$ with the filling
slope $ \phi_{\star}(\lambda_Z)=\mu$.

Let $e$ be an edge of $\Gamma$ with vertices $x$ and $y$.  For each
$N \in \DD(M, \Gamma)$ the edge $e$ corresponds to an incompressible
torus $T_e \subset N$, and its two vertices to two JSJ pieces $X$
and $Y$ of $N$, adjacent to the torus $T_e$. Denote the component of
$\partial X$ (resp. $\partial Y$) corresponding to $T_e$ by
$\partial _e X$ (resp. $\partial_e Y$). The embedded torus $T_e$
splits $N$ into two integral homology solid tori.

We call a slope on $\partial_e X$ {\it longitudinal} if it is equal
to the longitudinal slope $\lambda_V$ of an integral homology solid
torus $V$  bounded by $T_e$ in some $N \in \DD(M, \Gamma)$ and
containing $X$ (See Figure 1).

\psfrag{X}[c]{$X$} \psfrag{V}[c]{$V$}
\psfrag{H}[c]{$\mu_V=\phi^{-1}(\lambda_W)$}
\psfrag{I}[c]{$\lambda_V$} \psfrag{J}[c]{$\partial_eX$}
\psfrag{Y}[c]{$Y$} \psfrag{W}[c]{$W$}
\psfrag{C}[c]{$\mu_W=\phi(\lambda_V)$} \psfrag{B}[c]{$\lambda_W$}
\psfrag{D}[c]{$\partial_eY$} \psfrag{F}[c]{$\phi$}
\psfrag{1}[c]{$x$} \psfrag{2}[c]{$e$} \psfrag{3}[c]{$y$}
\psfrag{G}[c]{$\Gamma$}
\begin{center}
    \includegraphics[scale=.9]{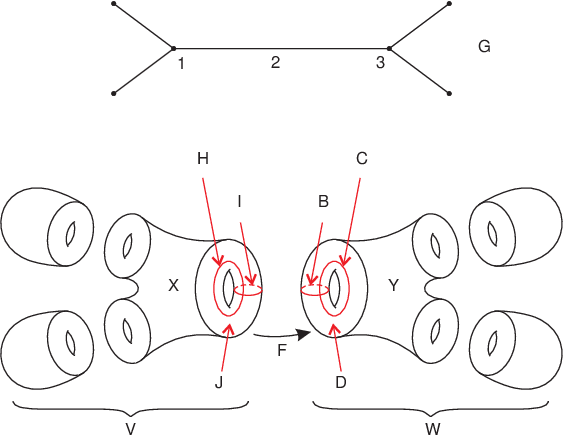}

    \centerline{Figure 1}
\end{center}

Let $N \in \DD(M, \Gamma)$. The incompressible torus $T_e$  splits
$N$ into two compact 3-manifolds $V$ and $W$ which are both integral
homology solid tori with boundary $T_e$: $N=V\cup_{T_e} W$. The fact
that $N=V\cup_\phi W$ is an integral homology sphere forces the
following:

\begin{Lemma}\label{gluing} The gluing map $\phi : \partial V \to \partial
W$ induces a map $\phi_{\star}$ on the first homology group, such
that $\phi_{\star}( \lambda_V)\cdot \lambda_W=\pm 1$ and
$\phi^{-1}_{\star}( \lambda_W)\cdot \lambda_V=\pm 1$.
\end{Lemma} \qed

\begin{Definition}

(1) We call a gluing map $\phi: \partial_e X \to \partial _e Y$
{\it allowable}, if there are two integral homology solid tori $V$ and $W$ such that  $(X, \partial_e X) \subset (V, \partial V)$, $(Y, \partial_e Y) \subset (W, \partial W)$,
and $N=V\cup_{ \phi} W \in \DD(M, \Gamma)$.

(2) An allowable gluing map $\phi: \partial_e X \to \partial _e Y$ is determined, up to isotopy, by the two pairs of slopes $\{\lambda_V, \mu_V \} \in H_1(\partial _e X; \Z) \times H_1(\partial _e X; \Z)$ and $\{\lambda_W, \mu_W\} \in H_1(\partial _e Y; \Z) \times H_1(\partial _e Y; \Z)$, such that:

\noindent (i)  $\lambda_V$ and $\lambda_W$ are the longitudinal slopes of the integral homology solid tori $V$ and $W$ which define the gluing map
$\phi$ to be allowable.

\noindent (ii) $\mu_V=\phi^{-1}_{\star}( \lambda_W)$ and $\mu_W=\phi_{\star}( \lambda_V)$. They are called longitudinal-images or $\ell$-images for short.

\noindent By Lemma \ref{gluing},  the pair $\{\lambda_V, \mu_V \}$
defines a basis of $ H_1(\partial _e X; \Z)$ and the pair
$\{\lambda_W, \mu_W\}$  a basis of  $ H_1(\partial _e Y; \Z) $. The
pairs of slopes $\{\lambda_V, \mu_V \}$ and $\{\lambda_W, \mu_W\}$
are called {\it gluing patterns} for the tori $\partial_e X$ and
$\partial_e Y$ (see Figure 1).

(3) Let $Y$ a vertex manifold of $\Gamma$ with $k$ boundary
components $\partial_i Y, i=1,...,k$. A {\it gluing pattern}  for
$\partial Y$ is a system of pairs of slopes $\{(\lambda_1,
\mu_1),..., (\lambda_k,\mu_k)\}  \in   (H_1(\partial _1 Y; \Z)
\times  H_1(\partial _1 Y; \Z)) \times \dots \times  (H_1(\partial
_k Y; \Z) \times  H_1(\partial _k Y; \Z))$ for which there are a
collection $Z_1, \dots , Z_k$ of integral homology solid tori and
gluing maps $\phi_i : \partial Z_i \to \partial_i Y, \, \, \, i= 1,
\dots , k$ such that (see the top picture of Figure 2):
\begin{enumerate}
\item $N=Y\cup_{\phi} \cup_{i=1}^k Z_i$ belongs to $\DD(M, \Gamma)$, with $\phi = \cup_{i=1}^k \phi_i : \cup_{i=1}^k \partial Z_i \to \partial Y$.
\item  $\lambda_i = \lambda_{W_i}$, where $W_i$ is the integral homology solid torus $N \setminus \text{int}( Z_i)$.
\item $\mu_i=\phi_{\star}( \lambda_{Z_i})$, for $ i= 1, \dots , k$ .

\end{enumerate}
\noindent Hence, each gluing map $\phi_i : \partial Z_i \to \partial_i Y$ is allowable and each pair of slopes $(\lambda_i, \mu_i)$ is a gluing pattern for the component $\partial_i Y$.

\noindent The slopes $\{\mu_1,..., \mu_k \}$ are called a system of
$\ell$-images for $\partial Y$.

(4) Two systems of slopes on $\partial Y $ are $A$-equivalent if
there is a homeomorphism $\tau : (Y, \partial Y) \to(Y, \partial Y)$
which is a product of Dehn twists along properly embedded essential
annuli in $Y$ and sends one set to the other. In the same way two
gluing patterns for $\partial Y$ are $A$-equivalent, if there is
such a homeomorphism of $(Y, \partial Y)$ sending one to the other.
\end{Definition}

\begin{Remark}\label{finite-finite}  Let $Y$ be a compact irreducible orientable 3-manifold with boundary an union of incompressible tori, and let $Aut(Y)$ be the mapping class group of $Y$. By  \cite{Joh} the subgroup $\mathcal{A}(Y) \subset Aut(Y)$ generated by Dehn twists along essential
tori and proper annuli in $Y$ is of finite index in $Aut(Y)$.
Therefore the finiteness, up to homeomorphisms of $Y$, of systems of
slopes on $\partial Y$ (or gluing patterns for $\partial Y$) is
equivalent to the finiteness of their $A$-equivalence classes, since
Dehn twists along essential tori in $Y$ do not affect the slopes on
$\partial Y$.
\end{Remark}

\psfrag{1}[c]{$z_{i+1}$} \psfrag{2}[c]{$z_{i-1}$}
\psfrag{3}[c]{$z_i$} \psfrag{A}[c]{$Y$} \psfrag{B}[c]{$Y(\mu_i)$}
\psfrag{C}[c]{$W_i$} \psfrag{D}[c]{$N$} \psfrag{E}[c]{$W_i(\mu_i)$}
\psfrag{F}[c]{$Y(\mu_1,\dots,\hat\mu_i,\dots,\mu_k)$}
\psfrag{G}[l]{deg 1 map} \psfrag{H}[c]{$\mu_i$}
\psfrag{I}[c]{$\lambda_{z_i}$}
\begin{center}
    \includegraphics[scale=.7]{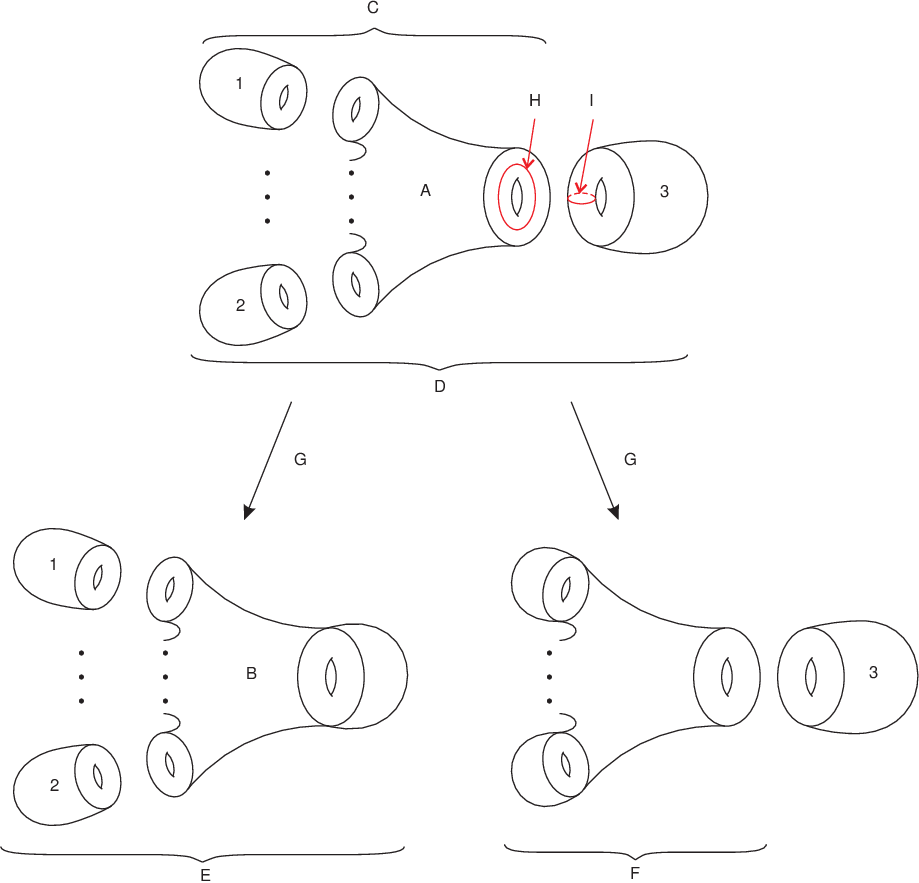}

    \centerline{Figure 2}
\end{center}

Now Proposition \ref{fixed graph} follows from the following result:

\begin{Proposition}\label{coord-finiteness} When $N$ runs over all elements in $\DD(M,
\Gamma)$, for each vertex manifold $Y$ of $\Gamma$ there are at most
finitely many $A$-equivalence classes of gluing patterns for
$\partial Y$, depending only on $M$.
\end{Proposition}

The first step of the proof of Proposition \ref{coord-finiteness} is given by
the following

\begin{Proposition}\label{finiteness of gluing slopes} When $N$ runs over all elements in $\DD(M,
\Gamma)$, for each vertex manifold $Y$ of $\Gamma$, there are at
most finitely many $A$-equivalence classes of $\ell$-images
$(\mu_1,..., \mu_k)$ on $\partial Y$, depending only on $M$.
\end{Proposition}

\demo Let $Y$ be a vertex manifold of $\Gamma$. By definition, for
each system of $\ell$-images $(\mu_1,..., \mu_k)$ on $\partial Y$,
there is a collection $Z_1, \dots , Z_k$ of integral homology solid
tori and gluing maps $\phi_i : \partial Z_i \to \partial_i Y, \, \,
\, i= 1, \dots , k$ such that: \begin{enumerate}
\item $N=Y\cup_{\phi} \cup_{i=1}^k Z_i$ belongs to $\DD(M, \Gamma)$, with $\phi = \cup_{i=1}^k \phi_i : \cup_{i=1}^k \partial Z_i \to \partial Y$.
\item $\mu_i=\phi_{\star}( \lambda_{Z_i})$, $i= 1, \dots , k$.
\end{enumerate}

Let $W_i$ be the integral homology solid torus $N \setminus
\text{int}(Z_i)$, then $Y$ is a JSJ-piece of $W_i$, for $i= 1, \dots
, k$.

We distinguish two cases according to whether $Y$ is hyperbolic or
Seifert fibered. By definition of $A$-equivalence, we prove the
finiteness of the systems of $\ell$-images in the hyperbolic case and the finiteness
of the systems of $\ell$-images, up to Dehn-twists along properly embedded essential annuli, in the  Seifer fibered case.

\noindent {\bf a) $Y$ is hyperbolic.} For each boundary component
$\partial_i Y$, $i\in \{1,..,k\}$, Thurston's hyperbolic Dehn
filling theorem \cite{Th1} shows that the manifold $Y( \mu_i)$
admits a complete  hyperbolic metric, except for a finite set of
slopes $\mu_i \in H_1(\partial _i Y; \Z)$, depending only on $Y$. So
we may assume that $Y(\mu_i)$ is hyperbolic.

Therefore $Y(\mu_i)$ is irreducible with incompressible boundary
tori, and it is a hyperbolic piece in the $JSJ$-decomposition of the
homology sphere $W_i(\mu_i)$. So $W_i(\mu_i)$ is an irreducible
homology sphere which is $1$-dominated by $N$ (see the top-left-down
picture of Figure 2), and thus dominated by $M$. Since $Y(\mu_i)$ is
a geometric piece of a manifold dominated by $M$, $Y(\mu_i)$ can
take only finitely many topological types, depending only on $M$ by
Theorem \ref{local domination}. Hence the hyperbolic volume of
$Y(\mu_i)$ takes finitely many values, depending only on $M$. Then
Thurston's hyperbolic Dehn filling theorem shows that $\mu_i$
belongs to a finite set of slopes in $H_1(\partial _i Y; \Z)$,
depending only on $Y$and $M$. Hence for each $i\in \{1,..,k\}$ there
are only finitely many possible $\ell$-images $\mu_i \in
H_1(\partial _i Y; \Z)$, depending only on $M$.

\noindent {\bf b) $Y$ is Seifert fibered.} Then the Seifert
fibration is unique, up to isotopy.

Suppose $Y=S(g, n-s;\frac{b_1}{a_1},\cdots,\frac{b_s}{a_s})$, where
$n-s=k$, and each $a_i>1, i=1,...,s$. Let $(c_i, h_i,) \in
H_1(\partial _i Y; \Z) \times  H_1(\partial _i Y; \Z)$ be a basis of
$H_1(\partial _i Y; \Z)$,
 where $h_i$ represents the fiber of the circle fibration induced on $\partial_i Y$ by the Seifert fibration of $Y$ and $c_i$ a section of this induced circle fibration. We set
$\mu_i=(a_{s+i}c_i + b_{s+i}h_i)$ in $H_1(\partial _i Y; \Z)$. There
is a degree one map from $N$ to the manifold $Y(\mu_1, ..., \mu_k)$
obtained by pinching each homology solid torus $Z_i$ to a solid
torus, hence  $Y(\mu_1, ..., \mu_k)$ is a Seifert fibered integral
homology 3-sphere. It follows that $g=0$, that is
$Y=S(n-s;\frac{b_1}{a_1},\cdots,\frac{b_s}{a_s})$. Moreover $n\ge 3$
since it is a JSJ piece of $N$. We further divide the discussion
into two cases:

Case (i):  $a_i\ne 0$ for each $i\in \{s+1, ..., s+k\}$. Then by
Lemma \ref{Seifert-Basic} (1) the Seifert fiberation of $Y$ extends
over the integer homology 3-sphere $Y(\mu_{1}, ..., \mu_k)$ and
$$\big(\prod_{i=1}^na_i \big)(\sum_1^n \frac{b_i}{a_i})=1.\qquad (*)$$

Like in case a), for each $i\in \{s+1, ..., s+k\}$, $W_i(\mu_i)$ is
an integer homology 3-sphere dominated by $M$. Moreover if
$a_{s+i}\ge 2$, the core of the filling solid torus becomes a
singular fiber of $Y( \mu_i)$ with index $a_{s+i}$, and thus
$Y(\mu_i)$ is a JSJ Seifert piece of $W_i(\mu_i)$. In this
case, $Y(\mu_i)$ can take only finitely many topological types,
depending only on $M$. Therefore $a_{s+i}\le C(M)$ for some integer
$C(M)$ depending only on $M$.

If $\{(a_{1}, b_{1}), ..., (a_{i}, b_{i}),...,(a_{j},
b_{j}),...,(a_{n}, b_{n})\}$ is a solution of the equation (*), then
for any integer $k$,  $\{(a_{1}, b_{1})$, ..., $(a_{i},
b_{i}+ka_i),...,(a_{j}, b_{j}-ka_j),...,(a_{n}, b_{n})\}$ is also a
solution of the equation (*). Those two solutions represent two systems of $\ell$-images on $\partial Y$
related by $k$ full Dehn twist (with sign) along an essential vertical annulus in $Y$
connecting $\partial_i Y$ and $\partial_j Y$, so they are in the same A-equivalence class.

We say that two solutions of (*) are in the same A-equivalence class, if one solution is obtained from another by finitely many Dehn twists along essential vertical annuli, like above.
It is an
elementary fact that there are only finitely many A-equivalence classes of
solutions for the equation (*) if each $a_i$ is bounded by a
constant $C(M)$. Hence  there are only finitely many A-equivalence classes of
systems of $\ell$-images on $\partial Y$, depending only on $M$.

Case (ii): $a_j=0$ for some $j\in \{s+1, ..., n\}$. Then by Lemma
\ref{Seifert-Basic} (2) $$b_j\prod_{i=1, i\ne j}^{ n} a_i=1 \qquad
(**).$$ This implies that $\mu_j=(0,1)$ and $\mu_i=(1,b_i)$ for
$i\ne j$. By performing $b_i$-full  Dehn twists along a vertical
annulus connecting $\partial_i Y$ and $\partial_j Y$ for each $i\ne j$,
we can transform the system of $\ell$-images $\{(1, b_1),..., (1, b_{j-1}),(0,1),(1,
b_{j+1}),...,(1, b_{n})\}$ to $\{(1, 0),..., (1, 0),(0,1),(1,
0),...,(1, 0)\}$. That is, for a fixed $j$, the system of $\ell$-images on $\partial Y$
 is unique, up to A-equivalence. Since $j$ picks only finitely many value, there are finitely many systems of $\ell$-images on $\partial Y$,
 up to  A-equivalence. \qed

The next step of the proof of Proposition \ref{coord-finiteness} is given by the following claim:

\begin{Claim}\label{claim:gluing pattern}
A gluing pattern $\{(\lambda_1, \mu_1),..., (\lambda_k,\mu_k)\}$ for
$\partial Y$ is uniquely determined by  the system of $\ell$-images
$(\mu_1,..., \mu_k)$ on $\partial Y$.
 \end{Claim}

\demo  Fix a system of $\ell$-images $(\mu_1,..., \mu_k)$ on $\partial Y$.
Recall that for this system of
slopes $(\mu_1,..., \mu_k)$, there is an homology sphere $N\in M(\Gamma)$
such that $N=Y\cup \cup_{i=1}^k Z_i$, and the image by gluing of each longitudinal
slope $\lambda_{Z_i}$ is $\mu_i$. As before, let $W_i=N\setminus Z_i$.

Then for each fixed $i\in \{1,...,k\}$, we can define a degree one
map $p: W_i\to Y(\mu_1,...,\hat \mu_i,...,\mu_k)$ which is the
identity on the boundary $\partial W_i=\partial _i Y$ by pinching
each $Z_j$ to a solid torus $U_j$ whose meridian is matched with
$\mu_j$, $j\in \{1,...,k,  j\ne i\}$ (see the top-right-down picture
of Figure 2). Then the Seifert surface $F_{W_i}$ is pinched to a
Seifert surface $F_*$ of the integral homology solid torus
$Y(\mu_1,...,\hat \mu_i,...,\mu_k)$ bounded by $\lambda_i$. Since
$Y(\mu_i,...,\hat \mu_i,...,\mu_k)$ is a fixed integral homology
solid torus, the longitudinal slope $\lambda_i \in
 H_1(\partial _i Y; \Z)$ is unique.
 This proves the claim.\qed

To finish the proof of  Proposition
\ref{coord-finiteness},  we distinguish two cases as usual:

\noindent {\bf Y is hyperbolic.} If $Y$ is hyperbolic, as we have
seen in the proof of Proposition \ref{finiteness of gluing slopes},
there are only finitely many systems of $\ell$-images $(\mu_1,...,
\mu_k)$ on $\partial Y$. So by Claim \ref{claim:gluing pattern}
there are only finitely many possible choices of gluing patterns
$\{(\lambda_1, \mu_1),..., (\lambda_k,\mu_k)\}$ for $\partial Y$.

\noindent {\bf  $Y$ is Seifert fibered}. Let $(\mu'_1,..., \mu'_k) =
\tau(\mu_1,..., \mu_k)$ be deduced from the system of $\ell$-images
$(\mu_1,..., \mu_k)$ by a homeomorphism $\tau: Y \to Y$ which is a
compositions of Dehn twists along vertical annuli. Then, by the
uniqueness in  Claim \ref{claim:gluing pattern}, the the system of
$\ell$-images $(\mu'_1,..., \mu'_k)$ determines the system
$(\lambda'_1,..., \lambda'_k) = \tau(\lambda_1,..., \lambda_k)$ of
longitudinal slopes on $\partial Y$. Hence an $A$-equivalent class
of systems of $\ell$-images on $\partial Y$ determine a unique
$A$-equivalent class of gluing pattern for $\partial Y$. Then by
Proposition \ref{finiteness of gluing slopes}, there are only
finitely many equivalent classes of gluing pattern for $\partial Y$.

This finishes the proof of  Proposition \ref{coord-finiteness}.\qed

The following corollary of Proposition \ref{coord-finiteness}
implies Proposition \ref{fixed graph}.

For each $N \in \DD(M, \Gamma)$, a submanifold $L\subset N$ is
called {\it canonical} if it is  a component of $N\setminus {\TT}$,
where ${\mathcal T}$ is a subfamily (may be empty) of  JSJ-tori of
$N$.

 \begin{Corollary}\label{cor:finiteness*} When $N$ runs over all elements in $\DD(M, \Gamma)$,
the canonical submanifolds of $N$ take at most finitely many
typological types, depending only on $M$.\end{Corollary}

 \demo The proof will be  by induction on the number $v(L)$ of
$JSJ$-pieces of a canonical submanifold $L$.

Corollary \ref{cor:finiteness*} is valid for $v(L)=1$ since
$\SS\HH(M)$ is finite. We suppose that it is
valid for $v(L)< m$ and  we are going to verify it for $v(L)=m$.

Fix a connected subtree
$\Gamma_*$ of $\Gamma$ with $m$ vertices  and let $D(M, \Gamma_*)$
be the set of canonical submanifolds with  dual $JSJ$-tree
$\Gamma_*$. Choose a vertex $y\in \Gamma_*$ with corresponding vertex
manifold $Y$.

For each canonical submanifold  $L\in D(M, \Gamma_*)$, we have
$L\setminus Y=\cup_{i=1}^{p} P_i$, where $P_i$ is a canonical
submanifolds and $v(P_i)<m$, hence $P_i$ can take only finitely many
topological types by the induction hypothesis.
 So we may fix the topology of $P_i$ for each $i = 1,..., p$.

We may suppose $\partial Y=\{\partial_1 Y,..., \partial_k Y\}$ and
$\partial_0 P_i$ is the component of $\partial P_i$ identified with
$\partial _i Y$ via a gluing map $\phi_i$ (reordering the components
of $\partial Y$ if needed). So we can rewrite
$L=Y\cup_{\{\phi_i\}_{i = 1,..., p}} \{P_i\}$.

 Fix a gluing pattern  $(\lambda_{P_i},\mu_{P_i}) \subset \partial
 _0 P_i$ for each $i = 1,..., p$. Then each gluing map  $\phi_i$ is determined  by the images
 $(\phi_i(\lambda_{P_i}), \phi_i(\mu_{P_i}))$ on $\partial _i Y$.
By definition $(\lambda_i= \phi_i(\mu_{P_i}), \mu_i= \phi_i(\lambda_{P_i}))$
 is a gluing pattern on $\partial_i Y$. Hence $\{(\lambda_{i},\mu_{i}), i=1,..., p\}$ forms a
 subset of a gluing pattern $\{(\lambda_{i},\mu_{i}), i=1,...,k\}$ for $\partial Y$. Any subset $\{(\lambda'_{i},\mu'_{i}), i=1,..., p\}$ of a $A$-equivalent gluing pattern
 $\{(\lambda'_{i},\mu'_{i}), i=1,...,k\}$ for $\partial Y$ provides a canonical submanifold $L'\in D(M,\Gamma_*)$ which is homeomorphic to $L$.
  By Proposition \ref{coord-finiteness}, there are
only finitely many $A$-equivalent classes of gluing patterns for $\partial Y$,
 depending only on $M$. Hence a canonical submanifold $L$ in $D(M,\Gamma_*)$ can take at most finitely many topological types, depending only on $M$.\qed

\section{Knot exteriors in $\S^3$}\label{knots}

By Theorem \ref{homology spheres} and an obvious twisted double
construction, one gets the following staigthforwards corollary:

\begin{Corollary}
Each compact orientable 3-manifold with a torus boundary 1-dominates
at most finitely many integral homology solid tori.

\end{Corollary}

A less direct and may be  more interesting result is the following

\begin{Theorem}\label{knot spaces}
A compact orientable 3-manifold $M$ dominates at most finitely many
exteriors of knots in $\S^3$.
\end{Theorem}

\demo We call the exterior $E(k)=S^3\setminus N(k)$ of a
knot $k$ in $\S^3$ a knot space, where $N(k)$ is a tubular
neighborhood of $k$ in $S^3$. The dual graph $\Gamma(k)$ to the
JSJ-decomposition of $E(k)$ is a rooted tree, where the root
corresponds to the unique vertex manifold containing $\partial
E(k)$.

Let $\KK(M)$ denote the set of homeomorphism classes of knot
spaces $E(k)$  dominated by $M$. By Lemma \ref{haken}, there are
only finitely many $\Gamma(k)$ for all $E(k)$ dominated by $M$. By
the local domination theorem (Corollary \ref{local domination with
boundary}) the JSJ-pieces of the knot spaces in $\KK(M)$ belong to a
finite set $\HH\SS(M)$.

For a given graph $\Gamma$, let $\KK(M, \Gamma) \subset \KK(M)$ be
the set of homeomorphism classes of  knot space $E(k)$, such
that:\begin{enumerate}

\item $E(k)$ is dominated by $M$.
\item The JSJ-graph $\Gamma(k)$ is abstractly isomorphic to
$\Gamma$.
\item Each vertex manifold has a fixed topological type.  Each torus
boundary component of the vertex manifold has assigned to an edge on
the vertex.
\end{enumerate}

Like in the case of integral homology speres, the proof of Theorem \ref{knot spaces} is reduced to the following:

\begin{Proposition}\label{fixed graph for knot}
The set $\KK(M, \Gamma)$ is finite.
\end{Proposition}

\demo  To apply the arguments in Section \ref{spheres} to the
present case, we will consider $S^3=E(k)\cup N(k)$ rather than just
consider $E(k)$. Precisely $\partial E(k)$ and the JSJ tori of
$E(k)$ provide an extended JSJ-splitting of $S^3=E(k)\cup N(k)$ with
one more torus $\partial E(k)$ and one additional  piece the solid torus $N(k)$.
The dual graph $\Gamma^*(k)$ of this extended decomposition is a
tree obtained by adding one leaf on the root.

Now for each JSJ-piece $Y$, different from $N(k)$,  of this extended JSJ decomposition of
$S^3=E(k)\cup N(k)$, we can define $A$-equivalent classes of gluing
patterns for $\partial Y \setminus \partial E(k)$ like in Section
\ref{spheres}. Similarly the proof of Proposition
\ref{fixed graph for knot} is reduced to the following

\begin{Proposition}\label{finiteness of gluing patterns-knot} When $N$ runs over all elements
in $\KK(M, \Gamma)$, for each vertex manifold $Y$ of $\Gamma$  there
are at most finitely many $A$-equivalent classes of gluing patterns
 for $\partial Y \setminus \partial E(k)$, depending only on $M$.
\end{Proposition}

And Proposition \ref{finiteness of gluing patterns-knot} is reduced
to the following

\begin{Proposition}\label{finiteness of gluning slopes-knot} When $N$ runs over all elements
in $\KK(M, \Gamma)$,  for each vertex manifold $Y$ of $\Gamma$ there
are at most finitely many $A$-equivalent classes of systems of
$\ell$-images
 on $\partial Y$, depending only on $M$.
\end{Proposition}

\begin{Remark} When $\partial E(k) \subset \partial Y$, in order to show the finiteness of $A$-equivalent classes of gluing patterns on $\partial Y \setminus \partial E(k)$,
we need  the finiteness of $A$-equivalent classes of systems of
$\ell$-images on $\partial Y$, including $\partial E(k)$.
\end{Remark}

\demo Let $Y$ be the given vertex manifold  with $k+1$ boundary
components $\partial_i Y, i=0,1,...,k$. Now $S^3\setminus Y=
\cup_{i=0}^{k} Z_i$, where $Z_0$ is a solid torus containing $N(k)$
and bounded by $\partial Y_0$, and $Z_i$ is a non-trivial knot space, bounded by $\partial
Y_i$ fore $i= 1,...,k$. Recall that for a system of $\ell$-images $(\mu_0, \mu_1,..., \mu_k)$ on
$(\partial _0 Y, \partial_1 Y, ..., \partial_k Y)$,   $\mu_0$ is the image of the boundary $\lambda_{Z_0}$ of a meridian disc of $Z_0$ and $\mu_i\subset
\partial_i Y$ is the image of the longitudinal slope $\lambda_{Z_i}$, $i=1,...,k$.

It is known that both the JSJ pieces in knot spaces and their gluing
are rather restrictive (see for example \cite[IX.22]{Ja},  \cite[
Chapter 2]{BS}). A Seifert JSJ-piece of a knot space is either a
torus knot space, a cable space or a composite space.

We may assume that $k \geq 1$, otherwise $Y$ is a hyperbolic or a
torus knot space and   by \cite{GL} the meridian $\mu_{0}$ is
unique. Below we distinguish three cases for the proof:

(i) {\bf Y is hyperbolic}. Since the the boundary tori $\partial Y
\setminus \partial_0 Y$ of the compact 3-manifold $Y(\mu_0)$ is
compressible, by \cite[Thm. 2.0.1]{CGLS}
 the $\ell$-image $\mu_0$ on $\partial_0 Y$ can belong to at most three distinct slopes. The argument for the finiteness of the remaining $\ell$-images $\mu_i$, $i=1,...,k$,
 is then the same as the corresponding part of the proof of Proposition \ref{finiteness of
gluing slopes}.

(ii) {\bf Y is a cable space}. Say $Y$ is a $(q,p)$-cable space with
$p \geq 2$. Then $Y$ is a Seifert fiber space over annulus with a
singular fiber of index $p$. Then $\partial Y = \partial_0 Y \cup \partial_1 Y$ and we choose a basis on $H_1(\partial Y_0; \mathbb Z)$ and $H_1(\partial Y_1; \mathbb Z)$ represented
 by a section of the circle fibrations induced on $\partial_0 Y$ and $\partial_1 Y$  by the Seifert fibration of $Y$ and the fiber of these induced circle fibrations.

 The fact that
$Y(\mu_0)$ is a solid torus forces $\mu_0$ to meet the fiber exactly
once, that is $\mu_0=(1, q_0)$ in $H_1(\partial Y_0; \mathbb Z)$.
Moreover $Y(\mu_1)$ must be a torus knot space $E(T_{q+sp, p})$,
which falls into $SH(M)$. Hence it has only finitely many
topological types. Therefore $\mu_1=(sp+q, q_1)$ and $s$ takes only
finitely many values. Then there are only finitely many
$A$-equivalent classes of systems of $\ell$-images on $\partial Y$
as in the corresponding part of the proof of Proposition
\ref{finiteness of gluing slopes}.

(iii) {\bf Y is a composite space}. It means that $Y$ is
homeomorphic to a product $S^1 \times D_k$ where $D_k$ is a disk
with $k$ holes. This corresponds to the case where the core $k_0$ of the solid torus $Z_0$ is
not a prime knot. In this case the $\ell$-image $\mu_{0}\subset
\partial_0 Y$ is isotopic to a fiber $h = S^1 \times \{\star \}$,
whose slope is determined by the topological type of $Y$. Then the
$A$-equivalence class of attaching patterns is unique as we shown in
the corresponding part of the proof of Proposition \ref{finiteness
of gluing slopes}.

This finishes the proof of Proposition \ref{finiteness of gluning
slopes-knot}. \qed

We call a homomorphism $\phi: \pi_1(M)\to \pi_1(N)$ between
3-manifold groups {\it non-degenerate}, if $\phi$ can be realized by
a proper map $f:M\to N$ of non-zero degree. The image of $\pi_1(M)$
by such a non-degenerate homomorphism has finite index in
$\pi_1(N)$.

Now we can translate Theorem \ref{knot spaces} into the following

\begin{Corollary}\label{knot group} The fundamental group of a compact, orientable 3-manifold
admits a non-degenerate homomorphism to only finitely many distinct
knot groups.
\end{Corollary}

Corollary \ref{knot group} is related to Simon's conjecture.

\begin{Conjecture}\label{epimorphism}  [Ki, Problem 1.12 (J. Simon)]
A knot group $\pi_1(S^3\setminus K)$ surjects onto at most finitely many distinct knot groups.\end{Conjecture}

This conjecture raised in 1970's has received recently a lot of
attention (see for example \cite{BBRW}, \cite{RW}, \cite{Si},
\cite{SW}, \cite{So4}). I. Agol and Y. Liu had confirmed Simon's
conjecture in the summer of 2010 \cite{AL} by proving that a
finitely generated group $G$ with the first Betti number
$\beta_1(G)=1$ surjects onto finitely many knot groups. Our result
holds with domain the fundamental group of any compact orientable
3-manifold and for non-surjective homomorphisms, but under the
restrictive condition that the homomorphism is non-degenerate.

We give now a criterion for a homomorphism between knot groups to be
non-degenerate:

\begin{Lemma}\label{lem:nondegenerate} A homomorphism $\phi :
\pi_{1}(E(k))\to \pi_{1}(E(k'))$ is non-degenerate iff it sends
the prefered longitude of $k$ to a non-trivial peripheral element of
$\pi_{1}(\partial E(k'))$.
\end{Lemma}

\demo On the boundary tori $\partial E(k)$ and $\partial E(k')$, let
$\{m, \ell \}$ and $\{m', \ell' \}$  be meridian-prefered longitude
pairs.

If $\phi$ can be realized by a proper map $f:E(k)\to E(k')$ of non-zero degree, then the restriction of $\phi: \pi_{1}(\partial E(k)) \to \pi_{1}(\partial E(k'))$ is injective, and thus
$\phi(\ell)$ is a non-trivial element in $\pi_{1}(\partial E(k'))$.

Conversely, assume $\phi(\ell)$ is a non-trivial element in
$\pi_{1}(\partial E(k'))$. It will be null-homologous in $H_1(E(k'))$,
hence $\phi(\ell) = \ell'^{n}$ with $n \in \Z \setminus \{0\}$.
Then $\phi(m)$ belongs to the centralizer of $\ell'^{n}$ in the
knot group $\pi_{1}(E(k'))$. By [JS, Chap. VI] and the description
of Seifert pieces in a knot complement, the centralizer of
$\ell'^{n}$ is the peripheral subgroup $\pi_{1}(\partial E(k'))$,
so $\phi(m)$ is a peripheral element which normally generates a finite index subgroup of the
knot group $\pi_{1}(E(k'))$, and so generates a finite index subgroup of its first homology group.
 Hence $\phi(m)$ must be equal to $pm' + q \ell'$ for some
integers $p \neq 0 , q \in \Z$. This shows that $\phi(\pi_{1}(\partial E(k))) \subset \pi_{1}(\partial E(k'))$
and that $\phi$ is injective on
$\pi_{1}(\partial E(k))$.

Then, since knot exteriors are $K(\pi,1)$-spaces, a standard
argument in algebraic topology and 3-manifold theory shows that the
homomorphism $\phi_i$  can be realized by a non-zero degree proper
map $f: E(k)\to E(k')$.  \qed

\begin{Remark} In  \cite{GR} \cite{HKMS} many examples of degenerate epimorphisms
between knot groups are given . There are epimorphisms between knot
groups which do not send a meridian to a meridian: Suppose a knot $k
\subset \S^3$ whose group $\pi_1(E(k))$ is normally generated by a
non-peripheral element $\mu$, see Lemma \ref{algebraic meridian}. By
\cite {Gon} there exists a knot $k' \subset \S^3$ and an epimorphism
from $\pi_1E((k'))$ onto $\pi_1(E(k))$ which sends a meridian of
$k'$ to $\mu$. The fact that knot groups are residually finite
\cite{He}, hence hopfian, and Property P for knots in $\S^3$
\cite{KM} imply that the knots $k$ and $k'$ must be distinct. This
construction has been pointed out to us by Cameron Gordon and Alan
Reid.
\end{Remark}

\begin{Lemma}\label{algebraic meridian}  Let $k$ be a $(1,1)$-knot in $\S^3$ which is not a $2$-bridge knot.
Then $\pi_1(E(k))$ is normally generated by a non-peripheral element.
\end{Lemma}

\demo Recall that a $(1,1)$-knot in $S^3$ is a knot which admits a
$1$-bridge presentation on a standard unknotted torus, therefore
$(1,1)$-knot is a tunnel number one knot and by construction the
fundamental group $\pi_1(E(k))$ is generated by two elements $a,m$
whith $m$  a meridian. Let $[a] = p[m] \in H_1(E(k); \mathbb Z)$,
then $\pi_1(E(k))$ is normally generated by the element $b =
am^{(1-p)}$. By \cite{BZ} this element cannot be peripheral since
$\pi_1(E(k))$ is generated by $b$ and $m$, and the fact that $k$ is
not a $2$-bridge knot. \qed

\medskip
Institut Math\'ematique de Toulouse, CNRS UMR 5219, Universit\'e
Paul Sabatier, 118 Route de Narbonne, F-31062 Toulouse Cedex 9,
France. E-mail address: boileau@math.univ-toulouse.fr

\medskip

Department of Mathematics and Statistics, University of Melbourne,
Parkville VIC. 3010, Australia. E-mail address:
rubin@ms.unimelb.edu.au

\medskip

LAMA Department of Mathematics, Peking University, Beijing 100871
China.\hfill \break \noindent E-mail address:
wangsc@math.pku.edu.cn

\end{document}